\documentclass[11pt]{amsart}

\usepackage{latexsym}
\usepackage{amssymb}
\usepackage{amscd}
\usepackage[all]{xy}
\usepackage[mathscr]{euscript}
\usepackage{dsfont}
\usepackage{hyperref}

\normalsize

\RequirePackage{xspace}

\newcommand{\thought}[1]{}
\renewcommand{\thought}[1]{ \textbf{[#1]}}

\usepackage{enumerate}        
\newenvironment{roenumerate}{\begin{enumerate}[\upshape (i)]}{\end{enumerate}}

\newcommand\nc {\newcommand}
\newcommand\rnc{\renewcommand}

\newcount\blopone
\newcount\xone
\newcount\xtwo
\newcount\ytwo

\newtheorem{theorem}{Theorem}[section]
\newtheorem{prop}[theorem]{Proposition}
\newtheorem{com}[theorem]{Comment}
\newtheorem{redu}[theorem]{Reduction}
\newtheorem{refinement}[theorem]{Refinement}
\newtheorem{summary}[theorem]{Summary}
\newtheorem{importnota}[theorem]{Important Notation}
\newtheorem{prblm}[theorem]{Problem}
\newtheorem{notation}[theorem]{Notation}
\newtheorem{defin}[theorem]{Definition}
\newtheorem{caution}[theorem]{Caution}
\newtheorem{remark}[theorem]{Remark}
\newtheorem{reminder}[theorem]{Reminder}
\newtheorem{illustration}[theorem]{Illustration}
\newtheorem{lemma}[theorem]{Lemma}
\newtheorem{construction}[theorem]{Construction}
\newtheorem{corollary}[theorem]{Corollary}
\newtheorem{example}[theorem]{Example}
\newtheorem{conclusion}[theorem]{Conclusion}
\newtheorem{triviality}[theorem]{Triviality}
\newtheorem{proto}[theorem]{Prototype Quasifibration}
\newtheorem{cauex}[theorem]{Cautionary Example}
\newtheorem{hypo}[theorem]{Hypothesis}
\newtheorem{subth}{ }[theorem]
\newtheorem{case}{Case}[theorem]
\newtheorem{ssubth}{ }[subth]
\newtheorem{subrmk}[subth]{Remark}
\newtheorem{facts}[theorem]{Facts}
\newtheorem{observation}[theorem]{Observation}

\nc\tri[1]{\begin{triviality}
\label{#1}}
\nc\obs[1]{\begin{observation}
\label{#1}
\begin{em}}
\nc\fac[1]{\begin{facts}
\label{#1}
\begin{em}}
\nc\cas[1]{\begin{case}
\label{#1}
\begin{em}}
\nc\rfn[1]{\begin{refinement}
\label{#1}}
\nc\prt[1]{\begin{proto}
\label{#1}}
\nc\lem[1]{\begin{lemma}
\label{#1}}
\nc\pro[1]{\begin{prop}
\label{#1}}
\nc\thm[1]{\begin{theorem}
\label{#1}}
\nc\cor[1]{\begin{corollary}
\label{#1}}
\nc\dfn[1]{\begin{defin}
\label{#1}}
\nc\sthm[1]{\begin{subth}
\label{#1}}
\nc\srmk[1]{\begin{subrmk}
\label{#1}
\begin{em}}
\nc\exm[1]{\begin{example}
\label{#1}
\begin{em}}
\nc\plm[1]{\begin{prblm}
\label{#1}
\begin{em}}
\nc\rmk[1]{\begin{remark}
\label{#1}
\begin{em}}
\nc\rmd[1]{\begin{reminder}
\label{#1}
\begin{em}}
\nc\ntn[1]{\begin{notation}
\label{#1}
\begin{em}}
\nc\smr[1]{\begin{summary}
\label{#1}
\begin{em}}
\nc\cau[1]{\begin{caution}
\label{#1}
\begin{em}}
\nc\hyp[1]{\begin{hypo}
\label{#1}}
\nc\imn[1]{\begin{importnota}
\label{#1}
\begin{em}}
\nc\rdn[1]{\begin{redu}
\label{#1}
\begin{em}}
\nc\cax[1]{\begin{cauex}
\label{#1}
\begin{em}}
\nc\cmt[1]{\begin{com}
\label{#1}
\begin{em}}
\nc\con[1]{\begin{construction}
\label{#1}
\begin{em}}
\nc\ill[1]{\begin{illustration}
\label{#1}
\begin{em}}
\nc\ssthm[1]{\begin{ssubth}
\label{#1}
\begin{em}}
\nc\cnc[1]{\begin{conclusion}
\label{#1}
\begin{em}}

\nc\elem{\end{lemma}}
\nc\erdn{\end{em}\end{redu}}
\nc\erfn{\end{refinement}}
\nc\eprt{\end{proto}}
\nc\ethm{\end{theorem}}
\nc\eobs{\end{em}\end{observation}}
\nc\ecor{\end{corollary}}
\nc\edfn{\end{defin}}
\nc\esthm{\end{subth}}
\nc\esrmk{\end{em}\end{subrmk}}
\nc\epro{\end{prop}}
\nc\etri{\end{triviality}}
\nc\eexm{\end{em}
\end{example}}
\nc\ecmt{\end{em}
\end{com}}
\nc\efac{\end{em}
\end{facts}}
\nc\ermk{\end{em}
\end{remark}}
\nc\ermd{\end{em}
\end{reminder}}
\nc\eill{\end{em}
\end{illustration}}
\nc\eplm{\end{em}
\end{prblm}}
\nc\ecas{\end{em}
\end{case}}
\nc\ecau{\end{em}
\end{caution}}
\nc\ecax{\end{em}
\end{cauex}}
\nc\eimn{\end{em}
\end{importnota}}
\nc\entn{\end{em}
\end{notation}}
\nc\econ{\end{em}
\end{construction}}
\nc\esmr{\end{em}
\end{summary}}
\nc\ehyp{
\end{hypo}}
\nc\ecnc{\end{em}
\end{conclusion}}
\nc\essthm{\end{em}
\end{ssubth}}

\nc\sst{\scriptstyle}
\newcommand{\comment}[1]{}
\newcommand{\ri}{\longrightarrow}
\newcommand{\sr}{\rightarrow}

\newcommand{\zz}{{\mathbb Z}}

\newcommand{\K}{{\mathbf K}}

\newcommand{\D}{{\mathbf D}}
\newcommand{\Dqc}{{\mathbf D_{\text{\bf qc}}}}

\newcommand{\pp}{{\mathbb P}}

\newcommand{\oo}{\otimes}

\nc\one{\mathds{1}}
\nc\op{^{\hbox{\rm\tiny op}}}
\nc\mth{^{\hbox{\rm\tiny th}}}

\nc\script{\mathscr}
\nc\z{\zeta}
\nc\bc{{\mathbb{BC}}}
\nc\ct{{\script T}}
\nc\cf{{\script F}}
\nc\cg{{\script G}}
\nc\cl{{\script L}}
\nc\cu{{\script U}}
\nc\cv{{\script V}}
\nc\ce{{\script E}}
\nc\ch{{\script H}}
\nc\cs{{\script S}}
\nc\car{{\script R}}
\nc\cd{{\script D}}
\nc\cc{{\script C}}
\nc\ck{{\script K}}
\nc\ca{{\script A}}
\nc\ci{{\script I}}
\nc\cj{{\script J}}
\nc\co{{\script O}}
\nc\cm{{\script M}}
\nc\cp{{\script P}}
\nc\cx{{\script X}}
\nc\cy{{\script Y}}
\nc\cz{{\script Z}}
\nc\bd{\begin{description}}
\nc\ed{\end{description}}
\nc\ctob{{\script C}at\big(\ci^{op},\ca\big)}
\nc\clim{{\ds\mathop{\rm lim}_{\ds\longleftarrow}}\,}
\nc\climi{\clim^{\!i}\,}
\nc\climn{\clim^{\!n}\,}
\nc\colim{{\ds\mathop{\rm colim}_{\ds\la}}}
\nc\oa{\overline{\ca}}
\nc\s{\sigma}
\nc\ta{\tau}
\nc\os{\overline\sigma}
\nc\ot{\overline\tau}
\nc\T{\Sigma}
\nc\Tm{\Sigma^{-1}}
\nc\de[1]{{\mathop{\rm deg(#1)}}}
\nc\Ad[1]{\mathop{\rm Ad}(#1)}
\nc\ad[1]{\mathop{\rm ad}(#1)}
\nc\kth{{\it K}--theory}
\nc\loc[1]{{\text{\rm Loc(#1)}}}
\nc\coloc[1]{{\text{\rm Coloc}(#1)}}

\def\der #1 {D\left(#1\right)}
\nc\prf{\begin{proof}}
\nc\eprf{\end{proof}}
\nc\ds{\displaystyle}
\nc\Tor{\text{\rm Tor}}

\nc\cb{{\script B}}
\nc\ab{{\script A}b}

\nc\be{\begin{roenumerate}}
\nc\ee{\end{roenumerate}}

\nc\cat[1]{{\script C}at\Big({\big\{#1\big\}}\op\,\,,\,\,\ab\Big)}
\nc\csab{{\script C}at\big(\cs^{op},\ab\big)}
\nc\ctab{{\script C}at\Big({\{\ct^\alpha\}}^{op},\ab\Big)}
\nc\csex{{\script E}x\big(\cs^{op},\ab\big)}
\nc\ctex{{\script E}x\Big({\{\ct^\alpha\}}^{op},\ab\Big)}
\nc\sub{\qquad\subset\qquad}
\nc\ctr[1]{{\left.\ct\left(-,#1\right)\right|}_{\cs}}
\nc\ctrf[2]{{\left.\ct\left(#1,#2\right)\right|}_{\cs}}
\nc\Ctr[1]{{\left.\ct\left(-,#1\right)\right|}_{\ct^\alpha}}
\nc\Ctrf[2]{{\left.\ct\left(#1,#2\right)\right|}_{\ct^\alpha}}

\nc\la{\longrightarrow}
\nc\nin{\noindent}
\nc\cad[1]{\text{card}(#1)}
\nc\eq{\quad=\quad}
\nc\BA{\begin{array}{c}}
\nc\EA{\end{array}}
\nc\barr{
\[
\begin{array}{cccccccccccccccc}
}
\nc\earr{
\end{array}
\]
}
\nc\as[1]{{\langle S\rangle}^{#1}}
\nc\sh{\text{\it shift}}

\nc\yy[1]{{\left.\ct\left(-,#1\right)\right|}_{\ct^c}}
\nc\vrep[2]{{\left.\ct\left(#1,#2\right)\right|}_{\ct^\alpha}}
\nc\da{\downarrow}
\nc\Hom{{\mathrm{Hom}}}
\nc\RHom{{\mathrm{RHom}}}
\nc\HHom{{\script H}{\mathrm{om}}}
\nc\RHHom{{\script{RH}}{\mathrm{om}}}
\nc\End{{\mathrm{End}}}
\nc\Ext{{\mathrm{Ext}}}
\nc\mr\Modtc
\nc\PExt{{\mathop{\rm PExt}}}
\nc\stm{\text{\rm stmod}(kG)}
\nc\stM{\text{\rm StMod}(kG)}
\nc\e{\varepsilon}
\nc\p{\mathfrak{p}}
\nc\q{\mathfrak{q}}

\nc\rs{\s^{-1}A}
\nc\br{{\{\s^{-1}A\}}}
\nc\ra\ri
\nc\y[1]{\mathbf{y}#1}
\nc\x[1]{\mathbf{z}#1}
\nc\mmod[1]{#1\text{--\rm mod}}
\nc\Mod[1]{#1\text{--\rm Mod}}
\nc\Md {\ensuremath{\mathop{\textup{Mod}}}}
\rnc\mod[1]{\ensuremath{\mathop{mod-\textup#1}}\xspace}
\nc\Modtc{\Mod{\ct^c}}
\nc\pgldim[1]{\mathop{\rm pgldim}\,#1}
\nc\tf{{\rm [TR5]}}
\nc\tfs{{\rm [TR5$^*$]}}
\nc\Fun{\text{\rm Funct}(F\op,\ab)}
\nc\sym{\text{\rm Sym}}
\nc\sgn{\text{\rm sgn}}
\nc\Pro{\text{\rm Prod}^{}_\alpha(F\op,\ab)}
\nc\Yt[1]{{\left.\Hom_\ct^{}\left(-,#1\right)\right|}_F^{}}
\nc\dl{\delta}
\nc\Proj[1]{#1\text{--\rm Proj}}
\nc\proj[1]{#1\text{--\rm proj}}
\nc\Flat[1]{\text{\rm Flat}\,#1}
\nc\Inj[1]{\text{\rm Inj\,}#1}
\nc\qc[1]{\text{\rm qc\,}#1}
\nc\ov{\overline}
\nc\wt{\widetilde}
\nc\wh{\widehat}
\nc\ph{\varphi}
\nc\tstr{{\it t}--structure}
\nc\spec[1]{{\text{\rm Spec}(#1)}}

\newcommand{\Dqcmi}{{\mathbf D^-_{\mathrm{qc}}}}
\nc\EProd{\text{\rm EProd}}
\nc\dfg{\mathbf{D}_{\mathrm{fg}}}
\nc\dcoh{\mathbf{D}^b_{\mathrm{coh}}}
\nc\dmcoh{\mathbf{D}^-_{\mathrm{coh}}}
\nc\ECoprod{\text{\rm ECoprod}}
\nc\Prod{\text{\rm Prod}}
\nc\Coprod{\mathrm{Coprod}}
\nc\COprod{\mathrm{coprod}}
\nc\ldimp{\text{\rm LDim}^{\prod}}
\nc\ldimc{\text{\rm LDim}^{\coprod}}

\nc\bn{\mathbf{n}}
\nc\bx{\mathbf{x}}
\nc\by{\mathbf{y}}
\nc\bz{\mathbf{z}}
\nc\bk{\mathbf{k}}

\nc\hoco{
\begin{picture}(40,10)
\put(20,0){\makebox(0,0)[b]{\text{\rm Hocolim}}}
\put(5,-2){\vector(1,0){30}}
\end{picture}\,}

\nc\holim{
\begin{picture}(40,10)
\put(20,0){\makebox(0,0)[b]{\text{\rm Holim}}}
\put(35,-2){\vector(-1,0){30}}
\end{picture}}

\nc\Cop{\text{\rm Coprod}}
\nc\seq{{\mathbb{S}_{\mathbf{e}}}}
\nc\se{{S^{\text{\tt e}}}}
\nc\te{{T^{\text{\tt e}}}}
\nc\LL{{\text{\bf L}}}
\nc\R{\text{\bf R}}
\nc\id{\text{\rm id}}
\nc\supp{\text{\rm supp}}
\nc\tsrt{\emph{t}--structure}
\nc\Loc{\mathrm{Loc}}
\nc\ovlp{\mathrm{overlap}}
\nc\codim{\mathrm{codim}}
\nc\add{\mathrm{add}}
\nc\Add{\mathrm{Add}}
\nc\Smr{\mathrm{smd}}
\nc\dperf[1]{\D^{\mathrm{perf}}(#1)}
\nc\Perf{\cp\text{\rm erf}}

\begin{document}

\author{Amnon Neeman}\thanks{The research was partly supported 
  by the Australian Research Council, and partly carried out while
  visiting the CRM in Barcelona. The author is grateful to both institutions
for their support}
\address{Centre for Mathematics and its Applications \\
        Mathematical Sciences Institute\\
        Building 145\\
        The Australian National University\\
        Canberra, ACT 2601\\
        AUSTRALIA}
\email{Amnon.Neeman@anu.edu.au}

\title{Strong generators in $\dperf X$ and $\dcoh(X)$}

\begin{abstract}
  We solve two open problems: first we prove a
  conjecture of Bondal and Van den Bergh, showing  
that the category $\dperf X$ is strongly generated whenever
$X$ is a quasicompact, separated scheme, admitting a cover by open affine
subsets $\spec{R_i}$ with each $R_i$ of finite global dimension.
We also prove that, for a noetherian scheme $X$ of finite
type over an excellent scheme of dimension $\leq2$, the derived
category $\dcoh(X)$ is strongly generated. The known results in this
direction all assumed equal characteristic, we have no such
restriction.

The method is 
interesting in other contexts: our key lemmas turn out to give a simple
proof that, if $f:X\la Y$ is a 
separated morphism
of quasicompact, quasiseparated
schemes such that $\R f_*:\Dqc(X)\la\Dqc(Y)$
takes perfect complexes to complexes of bounded-below
Tor-amplitude, then $f$ must be
of finite Tor-dimension.
\end{abstract}

\subjclass[2000]{Primary 18E30, secondary 18G20}

\keywords{Derived categories, schemes, compact generators}

\maketitle

\tableofcontents

\setcounter{section}{-1}

\section{Introduction}
\label{S0}

Let $\ct$ be a triangulated category and $G\in\ct$ an object. 
Bondal and Van den Bergh~\cite[2.2]{BondalvandenBergh04} made a string of
definitions,
regarding what it means for $G$ to generate $\ct$, 
and we briefly remind the reader.

\rmd{R0.-3}
Define the full subcategory $\langle G\rangle_1\subset\ct$ to consist
of  all direct summands of
finite coproducts of suspensions of $G$. For integers $n\geq 1$ we
inductively define subcategories $\langle G\rangle_{n}$: an object 
lies in $\langle G\rangle_{n+1}$ if it is a direct summand of an
object $y$ admitting a triangle $x\la y\la
z\la$ with
$x\in \langle G\rangle_1$ and $z\in \langle G\rangle_n$. The object
$G$ is a \emph{classical generator} if $\ct=\bigcup_{n=1}^\infty \langle
G\rangle_n$, and a 
\emph{strong generator} if there exists an $n$ with
$\ct=\langle
G\rangle_n$.
The category $\ct$ is called \emph{regular} if a strong
generator exists. Note that we are following the terminology
of Orlov~\cite{Orlov16}, elsewhere in the literature what we call \emph{regular} would
go by the name \emph{strongly generated.}
\ermd

\nin
Strong generators are particularly useful in triangulated categories
proper over a noetherian, commutative ring $R$.
We remind the reader of the definition of properness.

\rmd{R0.-1}
Let $R$ be a noetherian, commutative ring and let $\ct$ be an
$R$--linear triangulated category. The $R$--linearity of $\ct$ means that the Hom-sets
$\ct(X,Y)$ are all $R$--modules, and the composition maps
$\ct(X,Y)\times\ct(Y,Z)$ are all $R$--bilinear.

We say that the category $\ct$ is \emph{proper over $R$} if, for
each pair of objects $X,Y\in\ct$, the direct sum
$\oplus_{i=-\infty}^\infty\ct(X,\T^iY)$ is a finite $R$--module. This of
course implies that $\ct(X,\T^iY)$ vanishes for all but finitely many $i\in\zz$.
\ermd

A key theorem that shows the usefulness of the definitions above is:

\thm{T0.-5}
Let $R$ be a noetherian, commutative ring, let $\ct$ be a regular triangulated
category proper over $R$, and suppose
furthermore that $\ct$ is Karoubian, meaning idempotents split.
Then an $R$--linear functor $H:\ct\la\Mod R$ is representable if and only if
\be
\item
  $H$ is homological.
\item
  For any object $X\in\ct$,
  the direct sum $\oplus_{i=-\infty}^\infty H(\T^iX)$ is a finite
  $R$--module.
\ee
\ethm

\nin
 When $R$ is a field Theorem~\ref{T0.-5} is due to Bondal and Van den
Bergh~\cite[Theorem~1.3]{BondalvandenBergh04}, and
in the generality above it may be found in
Rouquier~\cite[Theorem~4.16 and Corollary~4.18]{Rouquier08}.

\exm{E0.-7}
In view of Theorem~\ref{T0.-5} it becomes interesting to find examples
of regular, Karoubian triangulated categories proper over a
noetherian ring $R$.
Let us begin with Karoubian: if $X$ is a quasicompact, quasiseparated
scheme then the category $\dperf X$, of perfect complexes over $X$,
is well-known to be Karoubian.

Now for examples of proper triangulated categories: if $X$ is of
finite type and separated over a noetherian ring $R$ we have the equivalence below,
which explains the terminology
\[
\{X\text{ proper over }R\}\Longleftrightarrow\{\dperf X\text{
  proper over }R\}.
\]
If $R$ is a field there is a simple, direct proof in
Orlov~\cite[Proposition~3.30]{Orlov16}.
The general case may be proved by similar
methods, using the curve-selection result in 
Lipman~\cite[Exercise~4.3.9]{Lipman09}. See~\cite[Corollary~4.3.2]{Lipman-Neeman07} or 
Lemma~\ref{L0.95} below for the argument.

For all the author knows there might be examples of
schemes $X$ for which $\dperf X$
is proper over $R$, without $X$ being of finite type and separated
over $R$.

Now it's time to discuss examples of regular triangulated
categories. This leads us to the first main
theorem of the article, conjectured by Bondal and
Van den Bergh (see Remark~\ref{R0.1.1.5}).
\eexm

\thm{T0.1}
Let $X$ be a quasicompact, separated
scheme. Then $\dperf X$ is regular
if and only if $X$ can be covered by open affine subschemes
$\spec{R_i}$ with each $R_i$ of finite global dimension.
\ethm

\rmk{R0.1.1}
In particular: if $X$ is a separated
scheme, of finite type over a noetherian ring
$R$ of finite Krull dimension, then
\[
\{X\text{ regular}\} \Longleftrightarrow\{\dperf X\text{ regular}\}.
\]
Combining Theorem~\ref{T0.1} with Example~\ref{E0.-7} tells us that,
as long as $X$ is a regular scheme proper over a noetherian ring $R$
of finite Krull dimension,
the category $\dperf X$ is Karoubian, regular and proper over $R$.
This gives lots of examples to which Theorem~\ref{T0.-5} applies.
\ermk

\rmk{R0.1.1.5}
We should say something about the history of Theorem~\ref{T0.1}.
Bondal and 
Van den Bergh~\cite[Theorem~3.1.4]{BondalvandenBergh04}
proved the special case where $X$ is smooth over a field $k$.
But
in the sentence immediately following
\cite[Theorem~3.1.4]{BondalvandenBergh04} they go on to say:
``Presumably the theorem is true under the weaker hypothesis that
$X$ is noetherian and regular''. Hence we view 
Theorem~\ref{T0.1}, in the generality proved here, as
having been conjectured by Bondal and 
Van den Bergh.

The case where $X$ is regular and of finite type over a field
$k$ follows from either 
Rouquier~\cite[Theorem~7.38]{Rouquier08}
or
Orlov~\cite[Theorem~3.27]{Orlov16}. Note that all three results
mentioned so far---that
is \cite[Theorem~3.1.4]{BondalvandenBergh04}, 
\cite[Theorem~7.38]{Rouquier08}
and
\cite[Theorem~3.27]{Orlov16}---assume
equal characteristic.

There is an old result by Kelly~\cite{Kelly65}, giving
Theorem~\ref{T0.1} in the special case 
where $X$ is affine.
Kelly's theorem is the only one prior to this article
which works in
mixed characteristic. But affine schemes are rarely proper over
a noetherian ring $R$, hence Kelly's theorem
does not produce interesting $\ct$'s to which
one could apply Theorem~\ref{T0.-5}. Our theorem is the first to
produce a large slew of geometric triangulated categories
which are Karoubian, regular, and proper over a noetherian
ring $R$ of mixed characteristic.
\ermk

The reader might be curious why one cares about results of the type
described above. One of the major consequences of such theorems is

\cor{C0.1.1.6}
Let $\ct$ be a regular, Karoubian triangulated category
proper over a noetherian, commutative ring $R$. Let $\cs\subset\ct$ be
a triangulated subcategory (in particular 
full). If the inclusion functor $I:\cs\la\ct$ has either a right or a
left adjoint then it has both adjoints, and $\cs$ is also Karoubian,
regular and 
proper over $R$. 
\ecor

For the purpose of the proof of Corollary~\ref{C0.1.1.6} we recall that,
if $A$ is a class of objects in $\ct$, then $A^\perp$ and $^\perp A$ are defined
to be the full subcategories of $\ct$ whose objects are
determined by the rules
\[
\begin{array}{ccc}
  A^\perp&=&\{t\in\ct\mid\ct(\T^ia,t)=0\,\,\forall i\in\zz,\,\forall a\in A\},\\
^\perp A&=&\{t\in\ct\mid\ct(t,\T^ia)=0\,\,\forall i\in\zz,\,\forall a\in A\}.
\end{array}
\]

\prf
The properties of being Karoubian, regular and proper over
$R$ are true for $\ct$ if and only if they are true for
$\ct\op$. Replacing $I:\cs\la\ct$ by $I\op:\cs\op\la\ct\op$ 
if necessary, we may assume that
$I$ has a right adjoint $I_\rho$. 

The fact that $\cs$ is proper over $R$ follows immediately from its
being a subcategory.
From the existence of the right adjoint to $I$ and  
\cite[Proposition~9.1.18 and Corollary~9.1.14]{Neeman99}
we have that $\cs={^\perp(\cs^\perp)}$. If
$e:s\la s$ is an idempotent in $\cs$ then it splits in the Karoubian
$\ct$, that is we have an isomorphism
$s\cong x\oplus y$ in $\ct$ taking $e$ to
$1_x^{}\oplus0_y^{}$. But $x\oplus y\cong s\in\cs={^\perp(\cs^\perp)}$
implies
that $x,y\in {^\perp(\cs^\perp)}=\cs$, that is $\cs$ is Karoubian. Since
\cite[Remark~9.1.15 and Theorem~9.1.16]{Neeman99}
give an equivalence of $\cs$ with the Verdier quotient 
$\cs\cong\ct/{\cs^\perp}$, if $G$ is a strong generator in $\ct$
then the image of $G$ is a strong generator in
$\ct/{\cs^\perp}\cong\cs$; therefore $\cs$ is also regular.

 Now let $t$ be an object of $\ct$. The functor
 $\ct\big(t,I(-)\big):\cs\la\Mod R$ is a homological functor, and
 because $\ct$ is proper over $R$ we have that, for any object
 $s\in\cs$, the direct sum $\oplus_{i=-\infty}^\infty
 \ct\big(t,I(\T^is)\big)$ is a finite $R$--module. By
 Theorem~\ref{T0.-5}, applied to the category $\cs$, this functor is
 representable: there exists an object $I_\lambda(t)\in\cs$ with
\[
\ct\big(t,I(s)\big)\cong \ct\big(I_\lambda(t),s\big)\ .
\]
The existence of the left adjoint to $I$ follows formally, it takes the object
$t\in\ct$ to $I_\lambda(t)\in\cs$.
\eprf

Corollary~\ref{C0.1.1.6} opens the door to the theory of
semiorthogonal decompositions with
admissible building blocks, we briefly remind the reader. If the
inclusion $I:\cs\la\ct$ has a right adjoint then
\cite[Remark~9.1.15]{Neeman99} says that the inclusion $J:\cs^\perp\la\ct$
has a left adjoint, giving $\ct$ a semiorthogonal 
decomposition $\ct=\langle\cs^\perp,\cs\rangle$. 
From Corollary~\ref{C0.1.1.6} we learn that 
$I$ and $J$ each have both a right and a left adjoint---the
subcategories $\cs$ and $\cs^\perp$ are both admissible.
We note that until now the
theory of semiorthogonal decompositions
in which the pieces are admissible has
been confined to equal characteristic, just because
all the interesting known examples, of $\ct$'s which are
Karoubian, proper over $R$ and regular, were in equal
charateristic.
 The point we want to
make here is that Remark~\ref{R0.1.1} changes this---it gives a
plethora
of examples in mixed characteristic.

\rmk{R0.1.1.6.5}
This remark is about the known constructions of examples,
we are interested in producing
\be
\item
  Karoubian, regular triangulated categories proper over a commutative,
noetherian ring $R$.
\setcounter{enumiv}{\value{enumi}}
\ee
In this remark we will assume the reader
is familiar with DG enhancements of triangulated categories---a reader
who doesn't feel comfortable with the theory can safely skip
ahead to Remark~\ref{R0.1.1.7},
we will never refer back to Remark~\ref{R0.1.1.6.5}.

So far we know two ways to obtain examples of (i):
\be
\setcounter{enumi}{\value{enumiv}}
\item
  If $R$ is of finite Krull dimension and $X$ is a
  separated, regular scheme proper
  over $R$, then 
  $\ct=\dperf X$ is an example of (i).
\item
  If $\ct$ is an example of (i) then so is any admissible subcategory
  $\cs\subset\ct$, that is any triangulated subcategory $\cs$ where
  the inclusion has a right or a left adjoint
  (and therefore both by Corollary~\ref{C0.1.1.6}).
\setcounter{enumiv}{\value{enumi}}
\ee
There is a classical construction we haven't mentioned so far,
which is based on the old
theorem of Kelly~\cite{Kelly65}.
\be
\setcounter{enumi}{\value{enumiv}}
\item
  If $\Lambda$ is a finite (possibly noncommutative) algebra over
  $R$ and $\Lambda$ is of finite global dimension, then the derived
  category $\D^b(\proj\Lambda)$ is an example of (i).
\setcounter{enumiv}{\value{enumi}}
\ee
As we've already noted, if the pair $\cs\subset\ct$ is as in (iii) then so
is $\cs^\perp\subset\ct$.
If $\ct$ is an example of (i) then so are $\cs$ and $\cs^\perp$. Of course
in this situation $\ct$ is a recollement of $\cs$ and $\cs^\perp$, see
\cite[Section~9.2]{Neeman99}. It is natural to wonder when this
process can be reversed. This means: suppose we have triangulated
categories $\cs$ and $\cs'$, both examples of (i). How can one
glue them to form a new category $\ct$, also an example of (i),
with fully faithful embeddings $I:\cs\la\ct$ and $J:\cs'\la\ct$,
which have right and left adjoints, and so
that $J(\cs')=I(\cs)^\perp$ and $I(\cs)={^\perp J(\cs')}$? Of course there
is the dumb gluing, namely $\cs\times\cs'$. But there are
many $\cs\subset\ct$ where the ``extension'' does not split, and the
question becomes how to reconstruct $\ct$ out of $\cs$ and $\cs'$.

At the level of triangulated categories no one has any idea how to do this.
But if we assume both $\cs$ and $\cs'$
have the structure of DG enhancements over the
commutative, noetherian base ring $R$ then there is a gluing
procedure which we briefly recall.

Because $\cs$ and $\cs'$ have (strong) generators
$G\in\cs$ and $G'\in\cs'$, we can let $S$ and $S'$ be the
DG $R$--algebras
$S=\End(G)$ and $S'=\End(G')$, where
the endomorphisms are understood in the DG enhancements.
We obtain equivalences
$\cs=\ch^0(\Perf-S)$ and $\cs'=\ch^0(\Perf-S')$. Given any DG $S$-$S'$ bimodule
$M$ we can form the DG matrix algebra over $R$
\[
T\eq \left(
\begin{array}{cc}
  S & M\\
  0 & S'
  \end{array}
\right)
\]
Then the triangulated category $\ct=\ch^0(\Perf-T)$ is Karoubian, there are
natural  fully faithful functors
$I:\cs\la\ct$ (with a right adjoint) and $J:\cs'\la\ct$ (with a left
adjoint), we have $J(\cs')=I(\cs)^\perp$ and $I(\cs)={^\perp J(\cs')}$, and
$\ct$ is regular. If the homology of the bimodule $M$ is finite
over $R$ then it's also true that $\ct$ is proper over $R$.

The
proofs of
most of the statements may be found in
Kuznetsov and Lunts~\cite{Kuznetsov-Lunts15}. For the proof that $\ct$
is regular whenever $\cs$ and $\cs'$ are see
Orlov~\cite[Proposition~3.20]{Orlov16}.

The construction of Kuznetsov and Lunts
is very much in the spirit of noncommutative algebraic
geometry, where one thinks of noncommutative schemes as
triangulated categories
$\cs=\ch^0(\Perf-S)$
for DG algebras $S$. And, without the hypotheses of regularity and properness
of $\cs$, one does indeed expect noncommutative behaviour. Against
this background comes a lovely theorem of
Orlov~\cite[Theorem~4.15 and Corollary~5.4]{Orlov16}, which
asserts
\be
\setcounter{enumi}{\value{enumiv}}
\item
  We're still interested in examples of (i). There are the ones that come
  from (ii) and (iv). And then we can produce
  more examples by forming admissible
  subcategories as in (iii) and using the gluing procedure of
  Kuznetsov and Lunts~\cite{Kuznetsov-Lunts15} sketched above.

  If the underlying commutative ring $R$ is a perfect field then
  any example of (i), obtainable by a finite iteration of these steps,
  can also be produced much more directly: we
  can form it by applying (iii) once to an
  example as in (ii). In particular we can avoid the noncommutative
  gluing procedure of Kuznetsov and Lunts, and the finite-dimensional
  algebra examples of (iv) are all special cases of the algebro-geometric
  examples obtainable from (ii) and (iii).
\setcounter{enumiv}{\value{enumi}}
\ee
I recommend Orlov's paper highly. The proof of the main theorem requires
producing new schemes, and admissible subcategories of their
categories of perfect complexes. The argument is a spectacular
display of the power of semiorthogonal decompositions,
coupled with the
equivalences
between admissible subcategories on different schemes
produced in Orlov~\cite{Orlov03}.
And the paper is also beautifully written: it begins
with a gentle introduction and a survey of the theory, which a
nonexpert like myself
found very helpful. In fact the current paper was born when
I was trying to understand Orlov~\cite{Orlov16}.
\ermk

\rmk{R0.1.1.7}
Note the generality of Theorem~\ref{T0.1}: we don't even assume the
schemes noetherian. The reader might wonder why we bother with this
level of
abstraction.

If we're willing to assume $X$ quasi-projective, over a noetherian
ring $R$, then the proof simplifies substantially. In fact the general
case is proved by reducing to the quasi-projective situation. As it
turns out the passage from the quasi-projective case to the general
one is not substantially simplified by assuming the schemes
noetherian---if we don't like the projectivity hypothesis, then we might as
well go whole hog and prove the theorem in the generality above.
\ermk

\rmk{R0.1.2}
We should note that one direction in Theorem~\ref{T0.1} is easy,
we show
\[
\left\{\begin{array}{c}X\text{ admits a cover by }\spec{R_i}\\
\text{ with }R_i\text{ of finite global dimension}\end{array}\right\} \Longleftarrow\{\dperf X\text{ strongly generated}\}.
\]
In fact more is true: if $\dperf X$ is strongly generated and 
$U=\spec R$ 
is any open affine subset of $X$, then $R$ must be of finite global
dimension.

We see this as follows: as $U$ is an open subset of $X$, the main
theorem of Thomason and Trobaugh~\cite{ThomTro} tells us that the
restriction
functor $j^*:\dperf X\la \dperf U$ is the idempotent completion of a
Verdier quotient map. If $G\in\dperf X$ is a strong generator it
follows that $j^*G\in\dperf U$ is a strong generator. But
for $U=\spec R$ we deduce that $R$ must be of finite global 
dimension---see for example
Rouquier~\cite[Theorem~7.25]{Rouquier08}.
Hence it only remains to prove the direction $\Longrightarrow$. 
\ermk

\rmk{R0.1.3}
If $X$ is noetherian the local hypothesis of Theorem~\ref{T0.1} is equivalent 
to requiring $X$ to be regular and of finite Krull dimension. But 
there are examples of nonnoetherian schemes $X$ satifying the hypotheses
of the theorem.
One source of examples is absolutely flat rings---rings for which 
every module is flat. From 
Salles~\cite[Proposition~3, page~702]{Salles81} we know that
for flat modules the projective dimension and pure
projective dimension agree, and hence for absolutely flat rings 
the global dimension (the supremum over all modules of 
their projective dimension) is equal to
the pure global dimension (the supremum over all modules of their 
pure projective dimension). 
Now \cite[Theorem~2.2]{Kielpinski-Simson75} or
\cite[Theorems~7.47 or 11.21]{Jensen-Lenzing} tell us that for
rings of cardinality $\leq\aleph_n$ the pure global dimension is
$\leq n+1$. Hence an absolutely flat ring of cardinality $\leq\aleph_n$
has global dimension $\leq n+1$.

Concretely: if $k$ is the field of two elements and $R$ is the ring 
$R=k[x_1^{},x_2^{},x_3^{}\ldots]/(x_1^{2}-x_1^{},x_2^{2}-x_2^{},x_3^{2}-x_3^{}\ldots)$
then $R$ is a nonnoetherian
 absolutely flat ring of cardinality $\aleph_0$, and its
global dimension is $1$. Thus $X=\spec R$ is a nonnoetherian
 example where Theorem~\ref{T0.1}
applies.
\ermk

We've had an extensive discussion of Theorem~\ref{T0.1} and its significance,
and it's about time to move on to the other major results in
the article. For the next major theorem
we recall the notion of a regular alteration
of a scheme.

\rmd{R0.1} 
Let $X$ be a noetherian scheme. A \emph{regular alteration} of $X$ is a 
proper, surjective morphism
$f:Y\la X$, so that
\be
\item
$Y$ is regular and finite dimensional.
\item
There is a dense open set $U\subset X$ over which $f$ is finite.
\ee
In 
\cite{deJong96,deJong97} de~Jong proves, among other things, that any
scheme $X$, separated and 
of finite type over an excellent scheme $S$ of dimension $\leq2$,
admits a regular alteration. Since de Jong's papers Nayak~\cite{Nayak09}
showed that any scheme $X$, separated and 
essentially of finite type over $S$, admits 
a localizing immersion into a scheme $\ov X$ 
separated and of finite type 
(even proper) over $S$. By restricting a regular alteration
of $\ov X$ to $X\subset\ov X$ it immediately follows that any scheme $X$, 
separated and essentially of finite type over an excellent scheme $S$
of dimension $\leq2$, also
has a regular alteration.

Observe that, if $X$ is separated and 
essentially of finite type over 
an excellent scheme $S$ of dimension $\leq2$, then so is every closed
subscheme of $X$. Therefore all closed subschemes of $X$ admit regular
alterations. Hence any $X$ which is separated and 
essentially of finite type over 
a separated, excellent scheme $S$ of dimension $\leq2$ satisfies 
\ermd

\hyp{H0.3}
A scheme $X$ satisfies Hypothesis~\ref{H0.3} if it is noetherian,
separated, and every closed subscheme
of $X$ admits a regular alteration.
\ehyp

\nin
Now we are ready to state our second main result:

\thm{T0.7}
Let $X$ be a scheme satisfying Hypothesis~\ref{H0.3}. Then the
triangulated category $\dcoh(X)$, of bounded
complexes of coherent sheaves on $X$, is regular.
\ethm

We owe the reader a summary of what was known in this direction.
Of course when $X$ is regular and finite-dimensional then
the inclusion $\dperf X\la\dcoh(X)$ is an equivalence,
and Theorem~\ref{T0.1} tells us that the equivalent categories
$\dperf X\cong\dcoh(X)$ are regular.
By a clever extension and
refinement of the argument in Bondal and
Van den Bergh~\cite[Theorem~3.1.4]{BondalvandenBergh04},
Rouquier~\cite[Theorem~7.38]{Rouquier08}
showed that $\dcoh(X)$ is regular whenever $X$ is a 
separated scheme 
of finite type over a perfect field $k$. 
The preprint by
Keller and 
Van den Bergh~\cite[Proposition~5.1.2]{Keller-Vandenbergh08}
generalized to separated schemes of finite type over arbitrary fields,
but this Proposition disappeared in the passage to the published 
version~\cite{Keller-Vandenbergh11}.
The reader might also wish to look at Lunts~\cite[Theorem~6.3]{Lunts10}
for a different approach to the proof.
If we specialize the result of Rouquier, extended by Keller and 
Van den Bergh, to the case where $X=\spec R$ is an affine scheme,
we learn that $\D^b(\mmod R)$ is regular whenever $R$ is of finite
type over a field $k$. In recent years there has been interest 
among commutative algebraists to understand this better: 
the reader is referred
to 
Aihara and Takahashi~\cite{Aihara-Takahashi11},
Iyengar and Takahashi~\cite{Iyengar-Takahashi15}, 
and Bahlekeh, 
Hakimian, Salarian
and Takahashi~\cite{Bahlekeh-Hakimian-Salarian-Takahashi15}
for a sample of the literature. There is also a connection
with the concept of the radius of the (abelian) category of
modules over $R$; see Dao and Takahashi~\cite{Dao-Takahashi14,Dao-Takahashi15}
and Iyengar and Takahashi~\cite{Iyengar-Takahashi15}.
The union of the known results
seems to be that $\D^b(\mmod R)$ is regular if $R$ is an
equicharacteristic  excellent
local ring, or essentially
of finite type over a field---see~\cite[Corollary~7.2]{Iyengar-Takahashi15}.
In~\cite[Remark~7.3]{Iyengar-Takahashi15} it is observed
that there are examples of commutative, noetherian 
rings for which $\D^b(\mmod R)$ is \emph{not} regular\footnote{In view of 
Theorem~\ref{T0.7}, we conclude that a ring $R$ with
$\D^b(\mmod R)$ not regular must be such that 
not every closed subscheme of $\spec R$
admits a regular alteration. Of course examples of 
noetherian rings $R$ with 
$\spec R$
not admitting a regular alteration are known: if $R$ is a
noetherian integral
domain admitting a regular alteration $f:X\la\spec R$, then over 
the generic point
of $\spec R$ the map $f$ is finite and faithfully flat. Hence $f$ 
must be finite and
faithfully flat over some nonempty open subset of $\spec R$, 
which by faithfully
flat descent would have to be regular. Thus if $R$ is a
noetherian integral domain 
which is not J-0, then $\spec R$ is a noetherian scheme not 
admitting a regular alteration.}.

What's really new about Theorem~\ref{T0.7} is that, unlike the predecessors
recalled above,
it works in mixed characteristic.

We should say something about our proofs.
It turns out to be convenient to work not in $\dperf X$ and $\dcoh(X)$ but
in the larger category $\Dqc(X)$; that is we switch to the unbounded derived
category of cochain complexes of sheaves of $\co_X^{}$--modules
with quasicoherent cohomology. 
There are unbounded versions of Theorems~\ref{T0.1}
and~\ref{T0.7}. We don't yet have the notation
to state these theorems, they will come in \S\ref{S27}: we will give
two theorems about $\Dqc(X)$, Theorems~\ref{T27.1} and~\ref{T27.5},
and show
\[
\text{\rm Theorem~\ref{T27.1}}\Longrightarrow
  \text{\rm Theorem~\ref{T0.1}}\,,\qquad\qquad
\text{\rm Theorem~\ref{T27.5}}\Longrightarrow
\text{\rm Theorem~\ref{T0.7}}\,.
\]
and it will remain to prove Theorems~\ref{T27.1} and~\ref{T27.5}.
For now we note a couple of facts about  Theorems~\ref{T27.1} and~\ref{T27.5}.
\be
\item
The hypotheses on $X$ are identical in
Theorems~\ref{T27.1} and~\ref{T0.1} and identical in
Theorems~\ref{T27.5} and~\ref{T0.7}. But the conclusions are about
$\Dqc(X)$ instead of $\dperf X$ or $\dcoh(X)$.
\item
  If $X$ is regular and satisfies Hypothesis~\ref{H0.3} then
  Theorems~\ref{T27.1} and~\ref{T27.5}
  give the same conclusion for $\Dqc(X)$.
\ee
As the
reader may have guessed, from (i) and (ii) above, the proof
of Theorem~\ref{T27.5} will
use regular alterations to deduce the result from
Theorem~\ref{T27.1}---Theorem~\ref{T27.1} 
establishes Theorem~\ref{T27.5} in the regular case,
and alterations are about reducing to the regular case.
In this field several
experts have tried to use de Jong's theorem, there is a slew of results
that are known in characteristic zero using resolution of singularities and
conjectured in characteristic $p$. Not surprisingly the
experts have thought
of the idea of employing alterations to
prove these conjectures.
What's remarkable is that Theorem~\ref{T27.5} is the first
success story, and the only one so far. It is an example---in 
\cite[proof of Theorem~4 (special case)]{Neeman16}
there is a simple proof of Theorem~\ref{T27.5}, deducing it
from Theorem~\ref{T27.1} for schemes all of whose
closed subschemes have resolutions of singularities.
In other words if we could use resolutions of singularities the
proof simplifies substantially.
The proof given here combines
the use of alterations with two old theorems
of Thomason's. Thomason's theorems fall into the area that nowadays
goes by the
name ``support theory''---it's the combination of regular
alterations and support
theory that proves Theorem~\ref{T27.5}.

We have said something about the key ideas that go into deducing
Theorem~\ref{T27.5} from Theorem~\ref{T27.1}, and it remains to discuss
the proof of Theorem~\ref{T27.1} a tiny bit.
We begin with a general definition:
let $\ct$ be any
triangulated
category with coproducts and let $G$ be an object. For integers $A\leq B$ we
define 
the full subcategory 
${\ov{\langle G\rangle}}_1^{[A,B]}\subset \ct$ 
to consist of all direct summands of arbitrary coproducts of objects in
the set $\{\T^{-i}G,\,A\leq i\leq B\}$---in other words we allow arbitrary
coproducts but restrict the permitted suspensions, only
suspensions in a prescribed range can occur. The inductive definition
is as with $\langle G\rangle_n$, an object belongs to
${\ov{\langle G\rangle}}_{n+1}^{[A,B]}$ if it is a direct summand of an object
$y$ admitting a triangle $x\la y\la
z\la$ with
$x\in {\ov{\langle G\rangle}}_1^{[A,B]}$ and $z\in {\ov{\langle
G\rangle}}_n^{[A,B]}$. Our next theorem, which is 
a sharpening of \cite[Theorem~4.1]{Lipman-Neeman07}
and the first key step in the
proof of Theorem~\ref{T27.1}, says

\thm{T0.3}
Let $m\leq n$ be integers, let $X$ be a scheme of finite type over
a noetherian ring $R$, and let $G$ be a classical generator for
$\dperf X\subset
\Dqc(X)$. Then there exist integers $N,A\leq B$ so that every object
$F\in\Dqc(X)^{\leq n}$ admits a triangle $D\la E\la F\la$ with
$E\in {\ov{\langle
G\rangle}}_N^{[A,B]}$ and $D\in\Dqc(X)^{<m}$.
\ethm

\rmk{R0.5}
Note the generality: we are trying to prove Theorem~\ref{T27.1} which, 
in the noetherian case,
is an assertion about regular schemes.
And I have just told the reader that the first key step
is a theorem which doesn't assume any regularity.
In passing we note that it is entirely possible that some of the hypotheses
on the scheme $X$ in
Theorem~\ref{T0.3} are superfluous. For all the author knows the
assertion might generalize to all quasicompact, quasiseparated schemes.
\ermk

Theorem~\ref{T0.3} is the
first ingredient we need in the proof of Theorem~\ref{T27.1}. It
enters into the proof of our next major result.

\thm{T0.11}
Let $j:U\la X$ be an open immersion of quasicompact, separated schemes.
If $H\in\Dqc(U)$ is a perfect complex and $G\in\dperf X$ a
classical generator, then for any triple of
integers $n,a\leq b$ there exists a triple of integers
$N,A\leq B$ so that
$\R j_*\ov{\langle H\rangle}_n^{[a,b]}\subset
\ov{\langle G\rangle}_N^{[A,B]}$.
\ethm
 
In the body of the article we will proceed as follows:
in \S\ref{S27} we state Theorems~\ref{T27.1} and~\ref{T27.5}
and prove that Theorem~\ref{T27.1} implies
Theorem~\ref{T0.1} and Theorem~\ref{T27.5} implies Theorem~\ref{T0.7}.
In \S\ref{S100} we show how to deduce Theorem~\ref{T27.5}
from Theorem~\ref{T27.1}, using regular alterations and support theory. Then
\S\ref{S2} and \S\ref{S3}
give the 
proof of Theorem~\ref{T0.3}---this is probably the hardest part, because we
don't want to make the simplifying assumption that $X$ is quasiprojective.
In \S\ref{S4} we wrap things up---first
we will prove Theorem~\ref{T0.11},
and finally we'll show
how to use Theorem~\ref{T0.11}, coupled with Kelly's old
theorem about regular
affine schemes, to finish off the proof of Theorem~\ref{T27.1}. It's
in the use of Kelly's theorem that regularity enters.

\ill{I0.13}
It's worth noting that Theorems~\ref{T0.3}
and~\ref{T0.11} are useful in other contexts, they have
applications having nothing to do with the regularity of
triangulated categories. We end the introduction with an
illustration.

We next give a simple proof
that a separated, finite-type morphism of noetherian schemes is perfect
if and only if it is proper and of finite Tor-dimension. The original
proof may be found in \cite[Theorem~1.2]{Lipman-Neeman07},
but the machinery
developed here makes the problem a triviality
and leads to sharper statements.
In passing I mention
that the reason I care, about different proofs of such statements,
is that I'd like to generalize to morphisms of stacks---at
the moment I don't know how to do this, I haven't yet proved
useful versions of Theorems~\ref{T0.3}
and~\ref{T0.11} valid for algebraic stacks.

Note that for simplicity
we assume the schemes noetherian; this assumption can be removed just
as in \cite[Theorem~1.2]{Lipman-Neeman07}, but without the noetherian
hypothesis the statements become
a little technical. In this paper we would rather not recall the definitions
of pseudocoherent and quasiproper morphisms.

We should perhaps remind the reader: a morphism of schemes $f:X\la Y$ is
called \emph{perfect} if $\R f_*:\Dqc(X)\la\Dqc(Y)$
takes perfect complexes to perfect complexes.
What is being asserted, in this Illustration, is that a 
finite-type, separated
morphism of noetherian schemes is perfect if and 
only if it is proper and of finite Tor-dimension.
One implication is classical, it's been known since
\cite[Corollaire~4.8.1]{Illusie71C}
that a proper morphism of finite Tor-dimension is perfect.
The converse is relatively recent. Anyway:
it's customary to break this up into two bits, dealing separately with
the properness and with the finite Tor-dimension.

If $f:X\la Y$ is proper, then
Grothendieck's old result~\cite[Th\'eor\`eme~3.2.1]{Grothendieck61} tells us that $\R f_*\dcoh(X)\subset\dcoh(Y)$.
And the sharp converse, which we prove below in Lemma~\ref{L0.95},
says that if $\R f_*\dperf X\subset\dcoh(Y)$ then $f$ must be
proper---you don't need to check what $\R f_*$ does to every
complex in $\dcoh(X)$, it's enough to check the perfect ones.
In Example~\ref{E0.-7} we've already mentioned a special
case of this lemma, where $Y$ is assumed
affine.  Lemma~\ref{L0.95} isn't new, it may be
found in~\cite[Corollary~4.3.2]{Lipman-Neeman07}, but the
nonnoetherian generality
of the statement and proof in~\cite{Lipman-Neeman07} might be confusing.

About finite Tor-dimension:
in Illusie~\cite[Corollaire~3.4.1(a)]{Illusie71C} we can find
the fact that, if $f:X\la Y$ is of finite
Tor-dimension, then $\R f_*$ takes
an object $H\in\Dqc(X)$ of bounded-below Tor-amplitude to an object
$\R f_*H\in\Dqc(Y)$ of bounded-below Tor-amplitude.
The sharp 
converse, see Proposition~\ref{P0.105} below, asserts that
if $\R f_*$ takes perfect complexes
to complexes of bounded-below Tor-amplitude then
$f$ is of finite Tor dimension. This result is new.
\eill

Now for the precise statements and proofs.

\lem{L0.95}
Let $f:X\la Y$ be a separated, finite-type map of noetherian schemes.
If $\R f_*$ takes $\dperf X\subset\Dqc(X)$ into $\dcoh(Y)\subset\Dqc(Y)$
then $f$ is proper.
\elem

\prf
By \cite[Exercise~4.3.9]{Lipman09} if suffices to show that
$\R f_*$ takes $\dcoh(X)$ to $\dcoh(Y)$.
Let $F$ be an object in $\dcoh(X)$, we wish to show that $\R f_*F$
belongs to $\dcoh(Y)$.
Shifting if necessary, we may assume $F\in\Dqc(X)^{\geq0}$.

Choose an integer $\ell$ so that $\R f_*\Dqc(X)^{\leq0}$ is
contained in $\Dqc(Y)^{\leq\ell}$.
By
\cite[Theorem~4.1]{Lipman-Neeman07}
we may choose a triangle $D\la E\la F\la$ with $E$ a perfect complex
and $D\in\Dqc(X)^{\leq-\ell-1}$. We deduce a triangle
$\R f_*D\la\R f_*E\la\R f_*F\la$ in $\Dqc(Y)$, and by the choice
of $\ell$ we know that $\R f_*D\in\Dqc(Y)^{\leq-1}$.
Hence the map $\R f_*E\la\R f_*F$ induces an
isomorphism of cohomology sheaves
in degrees $\geq0$.
The assumption
of the Lemma gives that $\R f_*E$ is in $\dcoh(Y)$, so
in degrees $i\geq0$ the cohomology sheaves
$\ch^i(\R f_*E)\cong\ch^i(\R f_*F)$ are coherent. But
$F$ belongs to $\Dqc(X)^{\geq0}$, hence
in degrees $i<0$ the cohomology sheaves $\ch^i(\R f_*F)$ vanish. Therefore
$\R f_*F$ belongs to $\dcoh(Y)$.
\eprf

\pro{P0.105}
Suppose $g:X\la Y$ is a separated morphism of
quasicompact, quasiseparated schemes.
If $\R g_*$ takes every perfect complex $F\in\Dqc(X)$ to an object
$\R f_*F\in\Dqc(Y)$
of bounded-below Tor-amplitude, meaning Tor-amplitude in the interval
$[\ell,\infty)$ for some integer $\ell$ which may depend on $F$,
  then $g$ must be of finite Tor-dimension.
\epro

\prf
First we note that the question is local in $Y$. Clearly it's local in
$Y$ to check that $g$ is of finite Tor-dimension, but the condition that
$\R g_*$ takes perfect complexes to complexes of
bounded-below Tor-amplitude does not at first
sight appear local. The hypothesis of the Proposition gives that
$\R g_*$ takes perfect complexes
to complexes of bounded-below Tor-amplitude. 
If we are given a cartesian square
\[
\xymatrix{
  \wt X \ar[r]^-{u} \ar[d]_-{f} & X\ar[d]^g \\
  \wt Y\ar[r]^-{v} & Y
}
\]
with $v$ an open immersion,
we need to show that
$\R f_*$ takes
perfect complexes to complexes of bounded below Tor-amplitude.

Therefore let $G$ be an
object in $\dperf{\wt X}$. Then $G\oplus \T G$ has a vanishing
image in $K_0(\wt X)$, and by the main theorem of
Thomason and Trobaugh~\cite{ThomTro} 
there exists an object
$H\in\dperf X$ with $u^*H\cong G\oplus \T G$. Base-change gives
an isomorphism $v^*\R g_*H\la\R f_*u^*H\cong \R f_{*}G\oplus\T\R f_{*}G$.
Since $\R g_*H$ is of bounded-below Tor-amplitude so is $v^*\R g_*H$, and
hence so is the direct summand $\R f_{*}G$.

We now know that the question is local in $Y$, hence we may assume $Y$ affine,
in particular separated.
We are given a separated morphism $g:X\la Y$; as $g$ and
$Y$ are separated so is $X$. Let $j:U\la X$ be an open immersion with
$U$ affine. Now apply Theorem~\ref{T0.11}
to the object $H=\co_U^{}\in\dperf U$ and any classical generator
$G\in\dperf X$---the existence of such a $G$ is
proved in \cite[Theorem~3.1.1]{BondalvandenBergh04} or
\cite[Theorem~4.2]{Lipman-Neeman07}. From
Theorem~\ref{T0.11} we have $\R j_*\co_U^{}\in
\ov{\langle G\rangle}_N^{[A,B]}$.
Hence
$\R g_*\R j_*\co_U^{}\in
\R g_*\ov{\langle G\rangle}_N^{[A,B]}\subset
\ov{\langle \R g_*G\rangle}_N^{[A,B]}$.
But by hypothesis the Tor-amplitude of $\R g_*G$ is bounded below, which
gives a uniform lower bound for Tor-amplitude of the objects
of $\ov{\langle \R g_*G\rangle}_N^{[A,B]}$, we only allow
suspensions in a range. Thus
$(gj)_*\co_U^{}=\R(gj)_*\co_U^{}\cong\R g_*\R j_*\co_U^{}$ is of bounded-below
Tor-amplitude, which means that at every point of the open affine
subset $U\subset X$ the map $g$ is of finite Tor-dimension. Since
$U$ is arbitrary $g$ is of finite Tor-dimension.
\eprf
\rmk{endofintro}
We should end the introduction with a list of obvious open problems
that follow from our results, and an account of the progress made on
these problems since the article was written. But first a word of
warning: to state some of the problems we assume familiarity with
the notation and lemmas that appear later in the article. Hence the reader
might wish to go on to the body of the article, before returning to
the open problems at the end.

We promised a list---here we go.
\be
\item
Is the hypothesis that $X$ is separated necessary in
Theorem~\ref{T0.-5} (the theorem characterizing those $X$ for which
$\dperf X$ is strongly generated), or in its unbounded 
version Theorem~\ref{T27.1}? 
Do the theorems 
generalize to all quasicompact, quasiseparated schemes?
We should note that much of the theory developed by Bondal and Van den Bergh
works in the generality of quasicompact, quasiseparated scheme, making the
question natural.

There has been some progress on this since this article
was written, see Jatoba~\cite{Jatoba19}.
\item
What about Theorem~\ref{T0.7} about the strong generation
of $\dcoh(X)$, or its unbounded version Theorem~\ref{T27.5}?
Is it necessary to assume that $X$ is separated, noetherian, and 
that every closed subscheme
has a regular alteration?

Here the beautiful article by Aoki~\cite{Aoki20} makes a giant step forward.
The scheme still has to be separated and noetherian, but it 
suffices for it to be
quasiexcellent and finite dimensional.  
The key idea is that, in place of combining de Jong's alterations 
with a couple of Thomason's old theorems as
in the current article (alterations
are not known to exist for all quasiexcellent schemes), Aoki
cleverly combines
Gabber's weak local uniformization theorem with the descendability
techniques of Mathew.
\item
Assume $X$ is a noetherian, quasiexcellent, finite dimensional 
separated scheme. By Aoki's theorem
the category $\dcoh(X)$ is regular. But $\dcoh(X)$ is the category of
compact objects in the compactly generated triangulated category
$\K(\Inj X)\cong\text{IndCoh}(X)$, see 
Krause~\cite[Proposition~2.3]{Krause05}
or Drinfeld and Gaitsgory~\cite[Proposition~3.4.2]{DrinfeldGaitsgory11}. 
Pursuing the
parallel with $\dperf X$ and the way its regularity was studied 
via the unbounded version Theorem~\ref{T27.1}, one can wonder whether
the obvious analog of Theorem~\ref{T27.1} holds for 
$\K(\Inj X)\cong\text{IndCoh}(X)$.
\item
What about the noncommutative versions, for example
the categories of modules over an Azumaya algebra? 
\item
One could look at fancier noncommutative versions, such as the 
noncommutative schemes of Kontsevich and Orlov. There again one
can ask about strong generation of $\dperf X$ and $\dcoh(X)$. Of
course $\dcoh(X)$ needs to be defined. As we will see in the remainder 
of this Remark, this has been 
done.
Which leads us to
\item
Are the techniques developed in the article useful in contexts 
unrelated to strong generation? For example to the problem of defining
$\dcoh(X)$ for a noncommutative scheme $X$, as in (v) above?
\ee
The answer to question (vi) is a resounding Yes, and the rest of
this Remark will be devoted to saying a tiny bit about the 
subject that has grown out of applying the techniques developed here.
For a fuller account of the theory the reader is referred to the survey
article~\cite{Neeman17B}, or the much more detailed research articles
\cite{Neeman17A,Burke-Neeman-Pauwels18,Neeman18,Neeman18A}. The
order in which the research articles are listed is logical,
they form a linear sequence of papers with each building on
its predecessors (including the current article).

Let $X$ be a scheme. In Definition~\ref{D3.3} we establish what it
means for a subcategory $\cs\subset\Dqcmi(X)$ to be 
\emph{approximable.} And then in the rest of Section~5 we produce
larger and larger approximable subcategories of the category
$\Dqcmi(X)$, culminating in Theorem~\ref{T3.9} which tells us that,
under mild conditions on $X$, the categories $\Dqc(X)^{\leq m}$
are approximable for every integer $m$,

The way this is generalized in the sequels is the obvious:
in place of $\Dqc(X)$ we allow any compactly generated 
triangulated category with a \tstr. More precisely: in the 
article~\cite{Neeman17A}, a compactly generated
triangulated category $\ct$ is defined to be approximable if $\ct$
possesses a \tstr\ such that the categories $\ct^{\leq m}$ satisfy
the obvious generalization of what holds in $\Dqc(X)^{\leq m}$ by
Theorem~\ref{T3.9}.

And then one proves many remarkable and useful structure theorems, which
hold for all
approximable triangulated categories. Among them is the
fact that the \tstr\ used in the definition of approximability is 
actually close to unique: any two must be equivalent
(under a suitable equivalence relation).
Not only do we have (up to equivalence) a canonical \tstr,
in any approximable triangulated
category $\ct$ one can define a string of intrinsic subcategories
$\ct^c$, $\ct^-$, $\ct^+$, $\ct^b$, $\ct_c^-$ and $\ct^c_b$, study
the relations among them, and prove useful theorems about the way 
they interact.

Since there is a survey of the results available in~\cite{Neeman17B}
let me stop soon---the reader wishing to go into more
detail can look at the survey. 
Here I only want to mention that, for $X$ a separated,
noetherian scheme and for $\ct=\Dqc(X)$, the intrinsic subcategory 
$\ct^b_c\subset\ct$
turns out to be $\dcoh(X)\subset\Dqc(X)$. And if $X$ is a noncommutative
scheme in the sense of Kontsevich or Orlov then, under mild
restrictions, the article 
\cite{Burke-Neeman-Pauwels18} proves that the category $\Dqc(X)$ 
is approximable. Thus starting with such a noncommutative
scheme we can declare $\dcoh(X)$ to be $\Dqc(X)^b_c$, and ask
when it is strongly generated. This is the sense in which I said,
in problem (v) above, that we know how to define $\dcoh(X)$ for
noncommutative $X$.
\ermk

\bigskip

\nin
{\bf Acknowlegements.}\ \ The author would like to thank 
Pieter Belmans, Daniel Bergh, Jack Hall, 
Dmitri Orlov, Olaf Schn\"urer and an anonymous 
referee
for many helpful suggestions that led to improvements on previous 
versions. Thanks also go to Oriol Ravet{\'o}s and 
Jan {\v{S}}{\v{t}}ov{\'{\i}}{\v{c}}ek for 
help with the relevant literature.

\section{Background}
\label{S1}

\rmd{R1.1}
Let $\ct$ be a triangulated category;
we begin by reminding the reader of some old definitions.
\be
\item
If $\ca$ and $\cb$ are two subcategories of $\ct$, then $\ca\star\cb$
is the full subcategory of all objects $y$ for which there exists a
triangle
$x\la y\la z\la$ with $x\in\ca$ and $z\in\cb$. 
\item
If $\ca$ is a subcategory of $\ct$, then $\add(\ca)$ is the full
subcategory containing all finite coproducts of objects in $\ca$. 
\item
If $\ca$ is a subcategory of $\ct$ and $\ct$ is closed
under coproducts, then $\Add(\ca)$ is the full
subcategory containing all (set-indexed) coproducts of objects in
$\ca$. 
\item
If $\ca$ is a full subcategory of $\ct$, then $\Smr(\ca)$ is the full
subcategory of
all direct summands of objects in $\ca$.
\ee
\ermd

\rmk{R1.3}
Reminder~\ref{R1.1}(i) is as in \cite[1.3.9]{BeiBerDel82}, while 
Reminder~\ref{R1.1}(iv) is identical with \cite[beginning of
2.2]{BondalvandenBergh04}.
Reminder~\ref{R1.1}(ii) and (iii) follow the usual conventions in
representation
theory; in \cite[beginning of
2.2]{BondalvandenBergh04} the authors adopt the
(unconventional) notation that $\add(\ca)$ and $\Add(\ca)$
are also closed under the
suspension---thus $\add(\ca)$ as defined in
\cite{BondalvandenBergh04} is what we would denote
$\add\big(\bigcup_{n=-\infty}^\infty\T^n\ca\big)$. 
 The definitions that follow are therefore slightly
different from \cite{BondalvandenBergh04}.

Recall that the octahedral lemma
gives $\ca\star(\cb\star\cc)=(\ca\star\cb)\star\cc$. The fact
that coproducts of triangles are triangles tells us that
\[
\left\{
\begin{array}{c}
\add(\ca)=\ca\\
\add(\cb)=\cb
\end{array}
\right\}\quad\Longrightarrow\quad
\{
\add(\ca\star\cb)=\ca\star\cb
\}\ .
\]
If $\ct$ is closed under coproducts then
\[
\left\{
\begin{array}{c}
\Add(\ca)=\ca\\
\Add(\cb)=\cb
\end{array}
\right\}\quad\Longrightarrow\quad
\{
\Add(\ca\star\cb)=\ca\star\cb
\}\ .
\]
Note that
the empty coproduct is $0$, hence $0\in\add(\ca)\subset\Add(\ca)$
for any $\ca$.
\ermk

\dfn{D1.5}
Let $\ct$ be a triangulated category and $\ca$ a subcategory. We
define subcategories 
\be
\item
$\COprod_1(\ca)=\add(\ca),\qquad\qquad\COprod_{n+1}(\ca)=\COprod_1(\ca)\star\COprod_n(\ca)$.
\item
$\Coprod_1(\ca)=\Add(\ca),\qquad\qquad\Coprod_{n+1}(\ca)=\Coprod_1(\ca)\star\Coprod_n(\ca)$.
\item
$\COprod(\ca)=\bigcup_{n=1}^\infty\COprod_n(\ca)$.
\item
$\Coprod(\ca)$ is the smallest strictly full subcategory of $\ct$
containing $\ca$ and
satisfying 
\[\Add\big(
\Coprod(\ca)
\big)\subset\Coprod(\ca),\qquad\qquad
\Coprod(\ca)\star \Coprod(\ca)\subset \Coprod(\ca).\]
\ee
\edfn

\rmk{R1.6.70707}
Definition~\ref{D1.5} is best viewed as
a useful technical refinement of Bondal and
Van den Bergh~\cite[the paragraphs right before and right
after Definition~2.2.1]{BondalvandenBergh04}. We do not allow direct
summands or arbitrary suspensions, but otherwise
$\COprod_{n}(\ca)$ is the analog of $\langle\ca\rangle_n$,
$\Coprod_{n}(\ca)$ is the analog of $\overline{\langle\ca\rangle}_n$,
$\COprod(\ca)$ is the analog of $\langle\ca\rangle$ and
$\Coprod(\ca)$ is the analog of $\overline{\langle\ca\rangle}$.
And consistency with 
Bondal and
Van den Bergh
also explains why Definition~\ref{D1.5}(iii) looks different
from Definition~\ref{D1.5}(iv).
\ermk

\obs{E1.7}
The definitions immediately give the inclusions
\[\xymatrix{
\COprod_n(\ca)  \ar@{^{(}->}[r] \ar@{^{(}->}[d] & \COprod(\ca) \ar@{^{(}->}[d] \\
\Coprod_n(\ca)  \ar@{^{(}->}[r] & \Coprod(\ca) 
}\]
In Remark~\ref{R1.3} we noted that $0\in\add(\ca)\subset\Add(\ca)$, which
we can rewrite as $0\in\COprod_1(\ca)\subset\Coprod_1(\ca)$. If $x\in\COprod_n(\ca)$ (respectively $x\in\Coprod_n(\ca)$),  the triangle
$0\la x\stackrel1\la x\la$ and Definitions~\ref{D1.5}(ii) and (iii) tell
us that $x\in\COprod_{n+1}(\ca)$ (respectively $x\in\Coprod_{n+1}(\ca)$).
That is
\[
\COprod_n(\ca)\subset\COprod_{n+1}(\ca),\qquad\qquad
\Coprod_n(\ca)\subset\Coprod_{n+1}(\ca).
\]
We have that
\[
\add\big(\COprod_1(\ca)\big)\eq\add\big(\add(\ca)\big)\eq\COprod_1(\ca),
\]
and similarly $\Add\big(\Coprod_1(\ca)\big)=\Coprod_1(\ca)$.
Induction on $n$ gives the equalities
\[
\add\big(
\COprod_n(\ca)
\big)=\COprod_n(\ca),\qquad\qquad\Add\big(
\Coprod_n(\ca)
\big)=\Coprod_n(\ca).
\]
Any finite set of objects of the increasing 
union $\COprod(\ca)=\bigcup_{n=1}^\infty\COprod_n(\ca)$ must lie in
$\COprod_n(\ca)$ for some large $n$, and the fact that 
$\add\big(
\COprod_n(\ca)
\big)=\COprod_n(\ca)$ tells us that the direct sum also lies in
$\COprod_n(\ca)$. Thus 
\[
\add\big(
\COprod(\ca)
\big)=\COprod(\ca).
\]
The associativity of the $\star$ operation gives 
\[
\begin{array}{rcl}
\COprod_m(\ca)\star\COprod_n(\ca)&=&\COprod_{m+n}(\ca),\\
\Coprod_m(\ca)\star\Coprod_n(\ca)&=&\Coprod_{m+n}(\ca),
\end{array}
\]
and the first of these identities tells us that 
$\COprod(\ca)=\bigcup_{n=1}^\infty\COprod_n(\ca)$ satisfies
\[
\COprod(\ca)\star\COprod(\ca)=\COprod(\ca).
\]
\eobs

We prove next the little lemma:

\lem{L1.9}
Let $\ct$ be a triangulated category and 
suppose we are given subcategories
$\ca,\cc,\cs,\cx,\cz\subset\ct$. Assume
$\add(\ca)=\ca$ 
and $\add(\cc)=\cc$. 
Suppose that
\be
\item
For any object $s\in\cs$, any morphisms
$s\la x$ and $s\la z$ with $x\in\cx$ and $z\in\cz$ factor as
\[
s\la a\la x,\qquad\qquad
s\la c\la z
\]
with $a\in \ca$ and $c\in \cc$.
\item
Any morphism $d\la x$, with $x\in\cx$ and
$d\in\big(\Tm\cc\big)*\cs$,
factors as $d\la a\la x$ with $a\in\ca$.
\ee
Then any morphism $f:s\la y$, with $s\in\cs$ and  
$y\in\cx\star\cz$, must factor as $s\la b\la y$ with $b\in \ca\star \cc$.
\elem

\prf
Because $y\in\cx\star\cz$ there exists a triangle $x\la y\la z\la\T x$ with
$x\in\cx$ and $z\in\cz$. The composite $s\la y\la z$ is a morphism from 
$s\in\cs$ to $z\in\cz$, and by (i) it
must factor as $s\la c\la z$ with $c\in\cc$. In 
other words we have a diagram where the row is a triangle and the square
commutes
\[\xymatrix@C+10pt{
 & s\ar[r] \ar[d]_f & c\ar[d] & \\
x\ar[r] & y\ar[r] & z\ar[r]& \T x
}\]
This we may complete to a morphism of triangles
\[\xymatrix@C+10pt{
d\ar[r]\ar[d]& s\ar[r] \ar[d]_f & c\ar[d]\ar[r] & \T d\ar[d]\\
x\ar[r] & y\ar[r] & z\ar[r]& \T x
}\]
By (ii) the morphism $d\la x$ 
factors as 
$d\la a\la x$ with $a\in\ca$. The composite $s\la c\la \T d$ vanishes
as these are two morphisms in a triangle---hence the longer 
composite $s\la c\la \T d\la\T a$ also vanishes.
 We obtain the following commutative diagram
\[\xymatrix@C+10pt{
 &  s\ar[dr]\ar@/^1pc/[drr]^0 \ar@/_1pc/[dd]_f & & \\
&  & c\ar[d]\ar[r] & \T a\ar[d]\\
x\ar[r] & y\ar[r] & z\ar[r]& \T x
}\]
The commutative square on the bottom right
may be extended to a morphism of triangles
\[\xymatrix@C+10pt{
a\ar[d]\ar[r]& \wt b\ar[d]_g\ar[r] & c\ar[d]\ar[r] & \T a\ar[d]\\
x\ar[r] & y\ar[r] & z\ar[r]& \T x
}\]
with $\wt b\in\ca\star\cc$, and the vanishing of $s\la c\la \T a$ tells us
that $s\la c$ must factor as $s\stackrel h\la \wt b\la c$. We have produced 
an object $\wt b\in\ca\star\cc$ and morphisms 
$s\stackrel h\la{\wt b}\stackrel g\la y$. A diagram chase shows that the two
composites $s\la z$ in the diagram
\[\xymatrix@C+10pt{
& s\ar[r]^h \ar[d]_f  & \wt b\ar[dl]^g &  & \\
x\ar[r] & y\ar[r] & z\ar[r]& \T x
}\]
are equal. It follows that $f-gh$ must factor as $s\la x\la y$. But $s\in\cs$ 
and $x\in\cx$, and (i) guarantees
that $s\la x$ factors as $s\la\wt a\la x$ with $\wt a\in\ca$.

The fact that $\add(\cc)=\cc$ implies that $0\in\cc$, and therefore
$\ca=\ca\star\{0\}\subset\ca\star\cc$. Therefore $\wt a\in\ca$ implies 
$\wt a\in\ca\star\cc$. Since we also have $\add(\ca)=\ca$ we deduce that
$\add\big(\ca\star\cc\big)=\ca\star\cc$, and the fact that $\wt a$ and $\wt b$
both lie in $\ca\star\cc$ means that so does $\wt a\oplus \wt b$. But now
the map $f-gh$ factors
through $\wt a$ and $gh$ factors through $\wt b$, hence
$f=(f-gh)+gh$ factors through $\wt a\oplus \wt b\in \ca\star\cc$.
\eprf

\rmk{R1.997381}
If $\cs\subset\ct$ is a triangulated subcategory and contains $\cc$, then
$\big(\Tm\cc\big)*\cs\subset\cs$ and 
hypothesis (ii) of Lemma~\ref{L1.9} follows from hypothesis (i). In this article
all the applications of Lemma~\ref{L1.9} will be in situations where $\cs$
is triangulated and contains $\cc$.

As it has turned out, in a sequel we will need to apply Lemma~\ref{L1.9}
in the generality of its statement---see \cite[Lemmas~2.2 and 2.3]{Neeman18}.
\ermk

\lem{L1.11}
Let $\ct$ be a triangulated category with coproducts, let $\ct^c$ be the 
subcategory of compact objects in $\ct$, and let $\cb\subset\ct^c$ be any
subcategory. Then
\be
\item
For $x\in\Coprod_n(\cb)$ and $s\in\ct^c$, any map $s\la x$ factors as
$s\la b\la x$ with $b\in\COprod_n(\cb)$.
\item
For $x\in\Coprod(\cb)$ and $s\in\ct^c$, any map $s\la x$ factors as
$s\la b\la x$ with $b\in\COprod(\cb)$.
\setcounter{enumiv}{\value{enumi}}
\ee
\elem

\prf
Let us first prove (i). We begin with the case $n=1$; any map $s\la x$,
with $s\in\ct^c$ and $x\in\Coprod_1(\cb)=\Add(\cb)$, is a map from the 
compact object $s$ to a coproduct of objects in $\cb$, and hence
factors through a finite subcoproduct. In particular it factors 
through an object of $\add(\cb)=\COprod_1(\cb)$.

Now suppose we know the theorem for all integers up to $n$. Apply 
Lemma~\ref{L1.9} with 
\[\cs=\ct^c, \quad\ca=\COprod_1(\cb),\quad\cc=\COprod_n(\cb),\quad\cx=\Coprod_1(\cb),\quad\cz=\Coprod_n(\cb).
\]
Induction tells us that the hypotheses of Lemma~\ref{L1.9} are satisfied,
hence any map from an object in $\ct^c=\cs$ to an
object in $\Coprod_{n+1}(\cb)=\cx\star\cz$ factors through
an object in $\COprod_{n+1}(\cb)=\ca\star\cc$.

It remains to prove (ii). Let $\car$ be the full subcategory of 
all objects $r\in\ct$ so that any map $s\la r$, with $s\in\ct^c$, factors
through an object in $\COprod(\cb)$. We need to show that 
$\Coprod(\cb)\subset\car$, and to do this we will prove three things. First:
\be
\setcounter{enumi}{\value{enumiv}}
\item
$\cb\subset\car$.
\setcounter{enumiv}{\value{enumi}}
\ee
This is obvious because any map $s\la b$, with $b\in\cb$, factors as
$s\la b\stackrel1\la b$. Next
\be
\setcounter{enumi}{\value{enumiv}}
\item
$\Add(\car)\subset\car$.
\setcounter{enumiv}{\value{enumi}}
\ee
\nin
\emph{Proof of (iv).}\ \ 
Suppose we are given an object $s\in\ct^c$, a set of objects 
$\{r_\lambda^{}\in\car,\,\lambda\in\Lambda\}$, and a morphism 
\[
\CD
\ds s @>f>> \ds\coprod_{\lambda\in\Lambda}r_\lambda^{}\ .
\endCD
\]
Because $s$ is compact $f$ can be factored as
\[
\CD
\ds s @>\Delta>>\ds\bigoplus_{i=1}^n s @>\oplus_{i=1}^n f_i>>
\ds\bigoplus_{i=1}^n r_i @>I>>\ds\coprod_{\lambda\in\Lambda}r_\lambda^{}
\endCD
\]
where $\Delta$ is the diagonal map and $I$ is the inclusion
of a finite subcoproduct. Because $s\in\ct^c$ and $r_i\in\car$, each map
$f_i:s\la r_i$ factors as $s\la c_i\la r_i$ with $c_i\in\COprod(\cb)$.
Hence the map $f$ factors as 
\[
\CD
\ds s @>>>
\ds\bigoplus_{i=1}^n c_i @>>>\ds\coprod_{\lambda\in\Lambda}r_\lambda^{}\ ,
\endCD
\]
and $\oplus_{i=1}^n c_i$ belongs to 
$\add\big(\COprod(\cb)\big)=\COprod(\cb)$.\hfill{$\Box$}
\be
\setcounter{enumi}{\value{enumiv}}
\item
$\car\star\car\subset\car$.
\setcounter{enumiv}{\value{enumi}}
\ee
\nin
\emph{Proof of (v).}\ \ 
Apply Lemma~\ref{L1.9} with $\cs=\ct^c$, $\ca=\cc=\COprod(\cb)$,
and $\cx=\cz=\car$. Any map $s\la y$, with $s\in\ct^c=\cs$ and
$y\in\car\star\car=\cx\star\cz$, must factor through an object in
$\ca\star\cc=\COprod(\cb)\star\COprod(\cb)=\COprod(\cb)$.
\hfill{$\Box$}

By definition $\Coprod(\cb)$ is the minimal subcategory of $\ct$ satisfying
(iii), (iv) and (v), hence $\Coprod(\cb)\subset\car$.
\eprf

\pro{P1.13}
Let $\ct$ be a triangulated category with coproducts, and let $\cb$ be a
subcategory of $\ct^c$. Then
\be
\item
Any compact object in $\Coprod_n(\cb)$ belongs to
$\Smr\big(\COprod_n(\cb)\big)$.
\item
Any compact object in $\Coprod(\cb)$ belongs to $\Smr\big(\COprod(\cb)\big)$.
\ee
\epro

\prf
Let $x$ be a compact object in $\Coprod_n(\cb)$ [respectively
 in $\Coprod(\cb)$]. 
The identity map 
$1:x\la x$ is a morphism from the compact object $x$ to $x\in \Coprod_n(\cb)$
[respectively to $x\in \Coprod(\cb)$],
and Lemma~\ref{L1.11}(i) [respectively Lemma~\ref{L1.11}(ii)] tells us 
that $1:x\la x$
factors through an object $b\in \COprod_n(\cb)$ [respectively 
$b\in \COprod(\cb)$]. Thus $x$ is a direct summand of $b\in 
\COprod_n(\cb)$
[respectively of $b\in\COprod(\cb)$].
\eprf

\lem{L1.15}
Let $\ct$ be a triangulated category with coproducts, and let $\cb$ be
an arbitrary subcategory. Then
\[
\Coprod_n(\cb)\sub\Smr\big(\Coprod_n(\cb)\big)\sub
\Coprod_{2n}(\cb\cup\T\cb)\ .
\]
\elem

\prf
The inclusion $\Coprod_n(\cb)\subset\Smr\big(\Coprod_n(\cb)\big)$
is obvious, we need to prove the inclusion
$\Smr\big(\Coprod_n(\cb)\big)\subset
\Coprod_{2n}(\cb\cup\T\cb)$. Assume therefore that $x$ is an object of
$\Smr\big(\Coprod_n(\cb)\big)$, that is there is an object $b\in
\Coprod_n(\cb)$ containing $x$ as a direct summand.
But then there is an idempotent map $e:b\la b$ whose image is $x$ and,
by \cite[proof of Proposition~1.6.8]{Neeman99},
$x$ is isomorphic to
the homotopy colimit of the sequence $b\stackrel e\la b\stackrel e\la 
b\stackrel e\la \cdots$. In other words there is a triangle
\[
\CD
\ds\coprod_{i=0}^\infty b @>>> x @>>> \ds\T\left[\coprod_{i=0}^\infty b\right] @>>>
\endCD
\]
where $\coprod_{i=0}^\infty b\in \Add\big(\Coprod_n(\cb)\big)=\Coprod_n(\cb)$,
and therefore $x$ belongs to
\begin{eqnarray*}
\Coprod_n(\cb)\star\Coprod_n(\T\cb)&\subset&
\Coprod_n(\cb\cup\T\cb)\star\Coprod_n(\cb\cup\T\cb)\\
&=&\Coprod_{2n}(\cb\cup\T\cb)\ .
\end{eqnarray*}
\eprf

\ntn{N1.17}
Let $\ct$ be a triangulated category with coproducts, and let $\cb\subset\ct$
be a subcategory. For any pair of integers
$m\leq n$ we will write
\[
\cb[m,n]=\bigcup_{i=-n}^{-m}\T^i\cb\ .
\]
In this notation, Lemma~\ref{L1.15} asserts that 
\[\Coprod_n(\cb)\sub\Smr\big(
\Coprod_n(\cb)
\big)\sub \Coprod_{2n}\big(\cb[-1,0]\big)\ .
\]
We permit $m=-\infty$ and/or $n=\infty$; for example $\cb[m,\infty)
=\bigcup_{i=-\infty}^{-m}\T^i\cb$.
\entn

In the Introduction, more precisely in the paragraph right
before Theorem~\ref{T0.3}, we introduced the 
subcategories  ${\ov{\langle
G\rangle}}_N^{[A,B]}$ and then went on
to express the results in terms of them. We are now in
a position to compare the
subcategories $ {\ov{\langle
G\rangle}}_N^{[A,B]}$ of the Introduction to
the subcategories $\Coprod_{N}\big(\cb[A,B]\big)$ of
this section. To do so, we adopt the notation that when
$\{G\}$ is the subcategory
with just one object $G$, we will write 
$\Coprod_{N}\big(G[A,B]\big)$ for what should more accurately be
denoted $\Coprod_{N}\big(\{G\}[A,B]\big)$.

\cor{C1.17.5}
For integers $N>0$, $A\leq B$ 
the identity 
${\ov{\langle
G\rangle}}_N^{[A,B]}
=\Smr\Big(
\Coprod_N\big(G[A,B]\big)
\Big)$
always holds. We furthermore have inclusions
\[\Coprod_N\big(G[A,B]\big)\sub{\ov{\langle
G\rangle}}_N^{[A,B]}
\sub \Coprod_{2N}\big(G[A-1,B]\big)\ .
\]
\ecor

\prf
When $N=1$ the identity ${\ov{\langle
G\rangle}}_1^{[A,B]}
=\Smr\Big(
\Coprod_1\big(G[A,B]\big)
\Big)$ is by the definitions of the two sides. For general $N$ the 
result follows by induction, since
\begin{eqnarray*}
{\ov{\langle
G\rangle}}_{N+1}^{[A,B]}&=&\Smr\Big[{\ov{\langle
G\rangle}}_{1}^{[A,B]}\star{\ov{\langle
G\rangle}}_{N}^{[A,B]}\Big]\\
&=&\Smr\Big[
\Smr\Big(
\Coprod_1\big(G[A,B]\big)
\Big)\star\Smr\Big(
\Coprod_N\big(G[A,B]\big)
\Big)
\Big]\\
&=& \Smr\Big[
\Coprod_1\big(G[A,B]\big)
\star
\Coprod_N\big(G[A,B]\big)
\Big]\\
&=&\Smr\Big[
\Coprod_{N+1}\big(G[A,B]\big)
\Big]\ ,
\end{eqnarray*}
where the third equality is by \cite[Lemma~2.2.1(ii)]{BondalvandenBergh04}.

It remains to prove the ``furthermore'' assertion of the Corollary.
Lemma~\ref{L1.15} asserts that, for any $\cb$, 
\[\Coprod_N(\cb)\sub\Smr\big(
\Coprod_N(\cb)
\big)\sub \Coprod_{2N}\big(\cb\cup\T\cb\big)\ ;
\]
applying this to $\cb=\{\T^iG,\,-B\leq i\leq-A\}$ we conclude
\[\Coprod_N\big(G[A,B]\big)\sub
{\ov{\langle
G\rangle}}_{N}^{[A,B]}
\sub \Coprod_{2N}\big(G[A-1,B]\big)\ .
\]
\eprf

\rmk{R1.17.7}
Corollary~\ref{C1.17.5} tells us that, for the purpose
of the statements
of Theorems~\ref{T0.3} and~\ref{T0.11} of
the Introduction, the subcategories ${\ov{\langle
G\rangle}}_{N}^{[A,B]}$ and  $\Coprod_N\big(G[A,B]\big)$ are 
interchangeable---either is contained in the other up to 
changing the integers $N,A,B$. As stated,  Theorems~\ref{T0.3}  
and~\ref{T0.11}
are assertions that there exist integers
$N,A,B$ so that the categories
${\ov{\langle
    G\rangle}}_{N}^{[A,B]}$
are large enough to contain certain objects.
Corollary~\ref{C1.17.5} implies that it 
is equivalent for there to
exist integers $N',A',B'$ so that
the categories $\Coprod_{N'}\big(G[A',B']\big)$ are large enough.

In the rest of the current article we will use only the categories
$\Coprod_N\big(G[A,B]\big)$, and never again mention ${\ov{\langle
G\rangle}}_{N}^{[A,B]}$.
We will prove Theorems~\ref{T0.3} 
and~\ref{T0.11} restated in terms of 
$\Coprod_N\big(G[A,B]\big)$, which happen to work better in iterating
the approximations that come up in the proofs.

The subcategories ${\ov{\langle
G\rangle}}_{N}^{[A,B]}$ are close to the historical
precursors
${\ov{\langle
G\rangle}}_{N}^{}$ considered by Bondal and Van den Bergh,
and the subcategories $\Coprod_N\big(G[A,B]\big)$ are a further 
refinement. As we have already said, back
in Remark~\ref{R1.6.70707},
we have introduced technical variants of
the definitions of Bondal and Van den Bergh that turn out
to be useful in certain proofs. But the truth is that most
of the time there is little need to consider the finest
version of these
invariants, and in subsequent articles we will revert to using
the coarser ${\ov{\langle
G\rangle}}_{N}^{[A,B]}$ in place of the finer $\Coprod_N\big(G[A,B]\big)$.
In the string of sequels 
\cite{Neeman17A,Burke-Neeman-Pauwels18,Neeman18,Neeman18A}
there is exactly one point at which a proof really makes use
the finer $\Coprod_N\big(G[A,B]\big)$,
with the cruder ${\ov{\langle
G\rangle}}_{N}^{[A,B]}$ being 
inadequate. See~\cite[Lemmas 4.3, 4.4 and 4.5]{Burke-Neeman-Pauwels18}.
\ermk

We have proved a number of general lemmas, and it is time to 
specialize a little. We end the section with an example.

\exm{E1.18}
Let $\ct=\D(\ca)$ for some abelian category $\ca$
satisfying AB4 (that is 
coproducts exist and are exact). Let $\cb\subset\ca$ be a subcategory,
which we view as embedded in $\D(\ca)$ in degree 0. If $m\leq n$ are integers
and $x^i,\,m\leq i\leq n$ are coproducts of objects
in $\cb$, then any cochain
complex of the form
\[\xymatrix{
\cdots\ar[r]& 0 \ar[r] &x^m_{}\ar[r] &
x^{m+1}_{}\ar[r] &\cdots\ar[r] &
x^{n-1}_{}\ar[r] & x^{n}_{}\ar[r] &0\ar[r] &\cdots
}\]
belongs to $\Coprod_{n-m+1}\big(\cb[m,n]\big)$. We see this by
induction on $n-m$: if $n-m=0$ the complex 
\[\xymatrix{
\cdots\ar[r]& 0 \ar[r] &x^m_{}\ar[r] &
0\ar[r] &\cdots
}\]
is a coproduct of objects in $\T^{-m}\cb$, that
is it belongs to $\Coprod_{1}\big(\cb[m,m]\big)$. The general case follows 
inductively,
by considering the triangles
\[\xymatrix{
\cdots\ar[r]& 0 \ar[r] &0\ar[d]\ar[r] &
0\ar[r]\ar[d] &\cdots\ar[r] &
0\ar[r]\ar[d] & x^{i}_{}\ar[r]\ar[d] &0\ar[r] &\cdots\\
\cdots\ar[r]& 0 \ar[r] &x^m_{}\ar[d]\ar[r] &
x^{m+1}_{}\ar[r]\ar[d] &\cdots\ar[r] &
x^{i-1}_{}\ar[r]\ar[d] & x^{i}_{}\ar[r]\ar[d] &0\ar[r] &\cdots\\
\cdots\ar[r]& 0 \ar[r] &x^m_{}\ar[d]\ar[r] &
x^{m+1}_{}\ar[r]\ar[d] &\cdots\ar[r] &
x^{i-1}_{}\ar[r]\ar[d] & 0\ar[r]\ar[d] &0\ar[r] &\cdots\\
 & & & & & &
}\]
where the top row belongs to $\Coprod_{1}\big(\cb[i,i]\big)$
while the bottom row belongs to $\Coprod_{i-m}\big(\cb[m,i-1]\big)$.
\eexm

\section{Unbounded versions of Theorems~\protect{\ref{T0.1}}
  and~\protect{\ref{T0.7}}}

\label{S27}

As mentioned in the Introduction, in this article we will provide
unbounded versions of Theorems~\ref{T0.1} and~\ref{T0.7}.
In this section we first state them, and then
prove that they imply Theorems~\ref{T0.1} and~\ref{T0.7}.
First for the unbounded version of Theorem~\ref{T0.1}:

\thm{T27.1}
Let $X$ be a quasicompact, separated
scheme. Suppose $X$ can be covered by open affine subschemes
$\spec{R_i}$ with each $R_i$ of finite global dimension. Then
there exists an object $G\in\dperf X$ and an integer $n>0$ with
$\Dqc(X)=\Coprod_{n}\big(G(-\infty,\infty)\big)$.
\ethm

\nin
The following result is well-known, but the proof is so simple that
we include it.

\lem{L27.3}
Theorem~\ref{T0.1} follows from Theorem~\ref{T27.1}.
\elem

\prf
Theorem~\ref{T0.1} is an if and only if statement, but in Remark~\ref{R0.1.2}
we noted that one direction is easy. It sufficed
to prove that, if $X$ satisfies the hypotheses of Theorem~\ref{T27.1}
and Theorem~\ref{T27.1} is known to be true,
then $\dperf X$ is regular.
Let $G\in\dperf X$ be the object whose
existence is guaranteed by Theorem~\ref{T27.1}.
Put $\cb=\{\T^iG,\,i\in\zz\}$. Then Theorem~\ref{T27.1}
tells us that $\Dqc(X)=\Coprod_n(\cb)$, and Proposition~\ref{P1.13}(i)
gives that $\dperf X=\Smr\big(\COprod_n(\cb)\big)$. This certainly implies
that $G$ strongly generates $\dperf X$.
\eprf

Now for the unbounded version of Theorem~\ref{T0.7}:

\thm{T27.5}
Let $X$ be a scheme satisfying Hypothesis~\ref{H0.3}. Then there
exists an object $G\in\dcoh(X)$ and an integer $n>0$ with
$\Dqc(X)=\Coprod_{n}\big(G(-\infty,\infty)\big)$.
\ethm

\nin
The fact that  Theorem~\ref{T0.7} follows from  Theorem~\ref{T27.5}
is not quite so immediate, we devote the rest of
the section to the proof.

\lem{L27.7}
Let $\ct$ be a triangulated category with coproducts and a 
\tstr\ so that $\ct^{\geq0}$ is closed under coproducts,
and let $G$ be
a bounded object of $\ct$. There is an integer
$M$ so that 
\[
\ct^{\geq n}\cap \Coprod_N\big(G(-\infty,\infty)\big)\sub
\Coprod_{2N}\big(G[n-M,\infty)\big)\ .
\]
\elem

\prf
Without loss of generality we may assume $n=0$, and replacing $G$ by
a suspension we may assume 
$G\in\ct^{\leq 0}\cap\ct^{\geq a}$ for some $a\leq 0$.
Let $\ci$ be the ideal of all maps $f:x\la y$ in $\ct$ 
so that any composite
$\T^iG\la x\stackrel f\la y$ vanishes, for any $i\in\zz$ and any
map $\T^iG\la x$. Now we observe
\begin{itemize}
\item
If $w\in\ct^{\geq m}$ is any object, there exists a triangle
$w'\la w\la w''\la$ in $\ct$ with $w\la w''$ in $\ci$,
with $w''\in\ct^{\geq m+a-1}$ and with $w'\in\Coprod_1\big(G[m,\infty)\big)$.
\end{itemize}
\nin\emph{Proof of $\bullet$.}\ \ 
We let $w'$ be the coproduct, over all nonzero maps $\T^iG\la w$, of $\T^iG$.
Let $w'\la w$ be the obvious map, and complete
to a triangle $w'\la w\la w''\la$. Since $w\in\ct^{\geq m}$ and 
$G\in\ct^{\leq 0}$ a nonzero map
$\T^iG\la w$ can only happen when $i\leq -m$, giving 
$w'\in\Coprod_1\big(G[m,\infty)\big)$, which also
means $w'\in\ct^{\geq m+a}$. The fact that $w\la w''$ belongs to
the ideal $\ci$ is immediate, and the triangle
$w\la w''\la\T w'$, coupled with the fact that
$w\in \ct^{\geq m}$ and $\T w'\in \ct^{\geq m+a-1}$, gives
$w''\in\ct^{\geq m+a-1}$ [recall that $a\leq0$].\hfill{$\Box$}

Next choose $x\in\ct^{\geq0}\cap\Coprod_N\big(G(-\infty,\infty)\big)$
and proceed inductively. Put $x_0^{}=x$, and let the map $x_i\la x_{i+1}$ be
the morphism $x_i\la x_i''$ of $\bullet$. Induction tells us that $x_i$
belongs to $\ct^{\geq i(a-1)}$, that the map $x_0\la x_{i+1}$
belongs to $\ci^{i+1}$, and that the object $w_{i+1}$ in the triangle
 $w_{i+1}\la 
x_0\la x_{i+1}$ belongs to $\Coprod_{i+1}\big(G[i(a-1),\infty)\big)$.
Since $x=x_0$ is assumed to belong to 
$\Coprod_N\big(G(-\infty,\infty)\big)$, an easy induction on $N$ 
shows that 
the map $x\la x_N^{}\in\ci^N$ must vanish. 
Hence the map $w_N^{}\la x$ is a 
split epimorphism, with $w_N^{}\in\Coprod_{N}\big(G[(N-1)(a-1),\infty)\big)$.
This makes $x$ a direct summand of an object
in $\Coprod_{N}\big(G[(N-1)(a-1),\infty)\big)$, and
Lemma~\ref{L1.15} now tells us that
$x$ belongs to  $\Coprod_{2N}\big(G[Na-N-a,\infty)\big)$.
\eprf

\lem{L27.9}
Let $X$ be a noetherian scheme and $G$ an object in $\dcoh(X)$. Any map
$f:E\la F$, with $E\in\dcoh(X)$ and $F\in\Coprod_N\big(G[M,\infty)\big)$,
factors through an object $F'\in\COprod_N\big(G[M,\infty)\big)$.
\elem

\prf
We prove this by induction on $N$, let us begin with the case
$N=1$. Since $G$ is bounded it belongs to some $\Dqc(X)^{\geq a}$, and
hence $\Coprod_1\big(G[M,\infty)\big)\subset\Dqc(X)^{\geq M+a}$. By
\cite[Theorem~4.1]{Lipman-Neeman07} there exists a triangle 
$C\la E\la K\la$ with $C$ compact and $K\in\Dqc(X)^{\leq M+a-1}$. The composite
$C\la E\stackrel f\la F$ is a map from a compact object $C$ to a coproduct 
$F=\coprod_{\lambda\in\Lambda}\T^{i_\lambda}G$, and hence factors through a
finite subcoproduct. That is we can write $F=F'\oplus F''$, with
$F'\in\COprod_1\big(G[M,\infty)\big)$ a finite subcoproduct, so that
the square in the diagram below commutes
\[\xymatrix@C+20pt@R+2pt{
C \ar[r]\ar[d] & E\ar[d]^f \ar[r] & K \ar[r] & \\
F'\ar[r] & F'\oplus F''\ar[r] & F''\ar[r] &
}\] 
The rows are triangles, hence we may complete to a morphism
of triangles 
\[\xymatrix@C+20pt@R+2pt{
C \ar[r]\ar[d] & E\ar[d]^f \ar[r] & K \ar[r]\ar[d] & \\
F'\ar[r] & F'\oplus F''\ar[r] & F''\ar[r] &
}\] 
But $K\in \Dqc(X)^{\leq M+a-1}$ and $F''$ is a summand of $F'\oplus F''\in
\Dqc(X)^{\geq M+a}$, and the map $K\la F''$ must vanish. Therefore the
map $E\la F=F'\oplus F''$ must factor through $F'
\in\COprod_1\big(G[M,\infty)\big)$.

We have proved the case $N=1$ of the Lemma, and
now it's time for the 
induction step.
Suppose therefore that we know the statement up to $N\geq1$. We apply
Lemma~\ref{L1.9} with $\ct=\Dqc(X)$, and with
 \[\cs=\dcoh(X),\quad
\ca=\COprod_1\big(G[M,\infty)\big),\quad
\cc=\COprod_N\big(G[M,\infty)\big),\]
\[
\cx=\Coprod_1\big(G[M,\infty)\big),\quad
\cz=\Coprod_N\big(G[M,\infty)\big).
\]
We deduce that any map $E\la F$, with $E\in\cs=\dcoh(X)$ and
$F\in\cx\star\cz=\Coprod_{N+1}\big(G[M,\infty)\big)$, must
factor through an $F'\in\ca\star\cc=\COprod_{N+1}\big(G[M,\infty)\big)$.
\eprf

\lem{L27.11}
Let $X$ be a noetherian scheme and 
let $G$ be an object in $\dcoh(X)$. Then
$\dcoh(X)\cap\Coprod_N\big(G(-\infty,\infty)\big)$ is contained
in 
\[
\Smr\Big[\COprod_{2N}\big(G(-\infty,\infty)\big)\Big].
\]
\elem

\prf
Observe that, for the standard \tstr\ 
on $\Dqc(X)$, we have 
$\dcoh(X)=\cup_{n=-\infty}^0\big[\dcoh(X)\cap\Dqc(X)^{\geq n}\big]$.
It therefore suffices to show that, for every $n\leq0$,
\[
\dcoh(X)\cap\Dqc(X)^{\geq n}
\cap\Coprod_N\big(G(-\infty,\infty)\big)
\sub\Smr\Big[\COprod_{2N}\big(G(-\infty,\infty)\big)\Big]\ .
\]
Lemma~\ref{L27.7}
gives us the inclusion
\begin{eqnarray*}
\dcoh(X)\cap\Dqc(X)^{\geq n}\cap\Coprod_N\big(G(-\infty,\infty)\big) 
&\subset&\dcoh(X)\cap\Coprod_{2N}\big(G[n-M,\infty)\big)
\end{eqnarray*}
for some integer $M$;
it suffices to show that any object
$x\in\dcoh(X)\cap\Coprod_{2N}\big(G[n-M,\infty)\big)$ 
belongs to $\Smr\Big[\COprod_{2N}\big(G(-\infty,\infty)\big)\Big]$.

To show this observe that the identity map $1:x\la x$
is  a map from $x\in\dcoh(X)$ to 
$x\in\Coprod_{2N}\big(G[n-M,\infty)\big)$, and by Lemma~\ref{L27.9}
it must factor through an object $b\in\COprod_{2N}\big(G[n-M,\infty)\big)$.
Thus $x$ is a direct summand of $b$, which belongs to
\[
\COprod_{2N}\big(G[n-M,\infty)\big)\sub
\COprod_{2N}\big(G(-\infty,\infty)\big)\ .
\]
\eprf

And now it's time to finish off.

\lem{L27.13}
Theorem~\ref{T0.7} follows from Theorem~\ref{T27.5}.
\elem

\prf
Suppose $X$ satisfies Hypothesis~\ref{H0.3}. Theorem~\ref{T27.5}
allows us to choose an object $G\in\dcoh(X)$ and an integer $N>0$
so that $\Dqc(X)=\Coprod_N\big(G(-\infty,\infty)\big)$.
 Therefore
\begin{eqnarray*}
\dcoh(X) &=& \dcoh(X)\cap \Dqc(X)\\
       &=& \dcoh(X)\cap\Coprod_N\big(G(-\infty,\infty)\big)\\
     &\subset& \Smr\big[\COprod_{2N}\big(G(-\infty,\infty)\big)\big]\\
   &\subset& \langle G\rangle_{2N}^{}\ ,
\end{eqnarray*}
where the first inclusion is by Lemma~\ref{L27.11} and the second inclusion
is obvious.
\eprf

\section{An object $G\in\dcoh(X)$ generating $\Dqc(X)$ in finitely many steps}
\label{S100}

This section is devoted to the proof that
Theorem~\ref{T27.5} follows from Theorem~\ref{T27.1}.
Let us begin with a general little lemma.

\lem{L1.1}
Let $(\ct,\oo,\one)$ be a tensor triangulated category and let
$H\in\ct$ be an object. The thick tensor ideal generated by $H$ is
the union, over all objects $C\in\ct$ and all integers $N>0$, of
the subcategories $\langle C\oo H\rangle_N^{}$.
\elem

\prf
Let $\ci$ be the thick tensor ideal generated by $H$. Because $\ci$ is
an ideal containing $H$
we have $C\oo H\in\ci$ for every
  $C\in\ct$, and because $\ci$ is thick $\langle C\oo H\rangle_N^{}\subset\ci$.
Therefore 
\[
\bigcup_{N>0,\,C\in\ct}\langle C\oo H\rangle_N^{}\sub\ci\ .
\]
We need to prove the reverse inclusion. For this it suffices to show
that $\cj=\bigcup_{N>0,\,C\in\ct}\langle C\oo H\rangle_N^{}$ is a
thick tensor ideal containing $H$.

Trivially $H\in\langle \one\oo H\rangle_1^{}\subset\cj$. 
Now let $K$ be an object in $\ct$. For
any object $C\in\ct$ and any 
integer $N>0$ we have 
\[
 K\oo\langle C\oo H\rangle_N^{}
\quad\subset\quad
\langle  K\oo C\oo H\rangle_N^{}\quad\subset\quad\cj\ ,
\]
and hence $K\oo\cj\subset\cj$.
Therefore $\ct\oo\cj\subset\cj$, that is $\cj$ is a tensor ideal.
The fact that $\T\langle C\oo H\rangle_M^{}=\langle C\oo H\rangle_M^{}$
tells us that $\T\cj=\cj$, and  
the inclusions
\begin{eqnarray*}
\langle C\oo H\rangle_M^{}\star \langle C'\oo H\rangle_N^{}
&\subset& \langle (C\oplus C')\oo H\rangle_M^{}\star \langle (C\oplus C')\oo H\rangle_N^{}\\
&\subset&\langle (C\oplus C')\oo H\rangle_{M+N}^{}
\end{eqnarray*}
imply that $\cj\star\cj\subset\cj$; it therefore follows
that is $\cj$ is triangulated.
Finally the fact that $\Smr\langle C\oo H\rangle_M^{}=\langle C\oo H\rangle_M^{}$
implies that $\Smr(\cj)=\cj$. Hence $\cj$ is a thick tensor ideal
containing $H$.
\eprf

And now we get down to business.

\medskip

\nin
\emph{Proof that Theorem~\ref{T27.5} follows from Theorem~\ref{T27.1}.}\ \ 
Let $X$ be a scheme satisfying
Hypothesis~\ref{H0.3}; in particular $X$ is noetherian. If $X$ does
not satisfy Theorem~\ref{T27.5}, the noetherian hypothesis allows us
to choose a closed subscheme $Z\subset X$ minimal among those which
do not satisfy Theorem~\ref{T27.5}. Replacing $X$ by $Z$, we may assume 
that all proper 
closed 
subschemes $Z\subset X$ satisfy Theorem~\ref{T27.5}---it suffices
to prove that so does $X$.

Next observe that we may assume $X$ reduced: let $j:X_{\mathrm{red}}\la X$
be the inclusion of the
reduced part of $X$, let $\ci$ be the ideal sheaf defining the 
subscheme $X_{\mathrm{red}}\subset X$, and let $n$ be an
integer so that $\ci^n=0$. Because $X$ is separated and noetherian 
\cite[6.7]{Bokstedt-Neeman93} 
tells us that the map $\D\big(\mathrm{QCoh}(X)\big)\la\Dqc(X)$ is an
equivalence---hence, up to replacing
an object $C\in\Dqc(X)$ by an isomorph, we may assume
any object $C\in\Dqc(X)$ is a complex of quasicoherent sheaves
on $X$. But then the complex $C$ admits a filtration
\[
0=\ci^nC\subset \ci^{n-1}C\subset\cdots\subset \ci C\subset C
\]
where the objects $\ci^jC/\ci^{j+1}C$ belong to $\R j_*\Dqc(X_{\mathrm{red}})$.
Thus $C$ belongs to 
\[
\big[\R j_*\Dqc(X_{\mathrm{red}})\big]^{\star n}=\underbrace{
\big[\R j_*\Dqc(X_{\mathrm{red}})\big]\star\big[\R j_*\Dqc(X_{\mathrm{red}})\big]
\star\cdots
\star\big[\R j_*\Dqc(X_{\mathrm{red}})\big]}_{n\text{ copies}}
\]
 and it suffices to 
prove that $\Dqc(X_{\mathrm{red}})=\Coprod_{\wt N}\big(\wt G(-\infty,\infty)\big)$ 
for some integer $\wt N>0$ and some object
$\wt G\in\dcoh(X_{\mathrm{red}})$.

Now let $f:Y\la X$ be a regular alteration.  Because $Y$ is 
finite-dimensional, separated and regular we may apply Theorem~\ref{T27.1}
to $Y$, after all we are assuming Theorem~\ref{T27.1}.
We may choose an object $G\in\dperf Y$ and an
integer $N$ so that $\Dqc(Y)=\Coprod_N\big(G(-\infty,\infty)\big)$.
Hence $\R f_*\Dqc(Y)=\R f_*\Coprod_N\big(G(-\infty,\infty)\big)\subset
\Coprod_N\big(\R f_*G(-\infty,\infty)\big)$ with $\R f_* G\in\dcoh(X)$.
The projection formula tells us that, for any object $C\in\Dqc(X)$, we
have $C\oo^\LL_X\R f_*\co_Y^{}\cong \R f_*\LL f^* C$, and hence
$\Dqc(X)\oo^\LL_X\R f_*\co_Y^{}\subset \Coprod_N\big(\R f_*G(-\infty,\infty)\big)$.

Now let us study the object $\R f_*\co_Y^{}$.
Since $X$ is reduced and the map 
$f$ is finite over a dense open set of $X$, there is a dense open subset
$V\subset X$ over which $f$ is finite and flat. Therefore the restriction
of $\R f_*\co_Y^{}$ to the open
set $V$ is a vector bundle---it is definitely a perfect complex.
By Thomason's localization theorem we may choose a perfect complex
$H$ on $X$, and a map $H\la \R f_*\co_Y^{}\oplus \T\R f_*\co_Y^{}$
inducing an isomorphism on $V$---see~\cite{ThomTro} or
\cite[statements~2.1.4 and 2.1.5]{Neeman92A}. 
Complete this map to a triangle
\[
\CD
H @>>> \R f_*\co_Y^{}\oplus \T\R f_*\co_Y^{} @>>> Q @>>>
\endCD
\]
Then $Q$ belongs to $\dcoh(X)$, and its restriction to $V$ vanishes.
By~\cite[Lemma~7.40]{Rouquier08}
there is an inclusion of a closed subscheme $i:Z\la X$, with
image contained in
$X-V$, and an object $P\in\dcoh(Z)$ so that $Q=\R i_*P$. Because
$Z$ is a proper closed subscheme of $X$ Theorem~\ref{T27.5} is true for $Z$,
and we may choose an object $G'\in\dcoh(Z)$ and an integer $M$ so
that $\Dqc(Z)=\Coprod_M\big(G'(-\infty,\infty)\big)$. Now let
$C\in\Dqc(X)$ be arbitrary and tensor the triangle above with $C$;
we obtain the triangle
\[
\CD
C\oo_X^\LL H @>>>C\oo_X^\LL\big[ \R f_*\co_Y^{}\oplus \T\R f_*\co_Y^{}\big] 
@>>> C\oo_X^\LL \R i_*P @>>>
\endCD
\]
In the previous paragraph we saw that 
$C\oo_X^\LL\big[ \R f_*\co_Y^{}\oplus \T\R f_*\co_Y^{}\big]$
belongs to the category $\Coprod_N\big(\R f_*G(-\infty,\infty)\big)$,
while $C\oo_X^\LL \R i_*P\cong \R i_*\big[\LL i^*C\oo_Z^\LL P\big]$
belongs to $\R i_*\Coprod_M\big(G'(-\infty,\infty)\big)\subset
\Coprod_M\big(\R i_*G'(-\infty,\infty)\big)$. The triangle tells
us that $C\oo_X^\LL H$ belongs to 
\[
\Coprod_M\big(\R i_*G'(-\infty,\infty)\big)\star 
\Coprod_N\big(\R f_*G(-\infty,\infty)\big)
\]
which is contained in 
\[
\Smr\Big[
\Coprod_M\big((\R f_*G\oplus \R i_*G')(-\infty,\infty)\big)\star 
\Coprod_N\big((\R f_*G\oplus \R i_*G')(-\infty,\infty)\big)
\big]
\]
By Lemma~\ref{L1.15} this is contained in 
$\Coprod_{2(M+N)}\big((\R f_*G\oplus\R i_*G')(-\infty,\infty)\big)$. Thus
\[
\Dqc(X)\oo_X^\LL H\subset 
\Coprod_{2(M+N)}\big((\R f_*G\oplus\R i_*G')(-\infty,\infty)\big)
\]
where $\R f_*G\oplus\R i_*G'$ belongs to $\dcoh(X)$.

Now let us study the object $H$. It is a compact object, and on the subset
$V$ it is isomorphic to the direct sum $\R f_*\co_Y^{}\oplus \T\R f_*\co_Y^{}$.
Since over $V$ the map $f$ is finite, flat and surjective, the
object $\R f_*\co_Y^{}$ restricts to a nowhere vanishing vector bundle 
on $V$. Thus the support of $H$ contains the dense open set $V$, 
and since 
the support of the compact object
$H$ is closed it must be
all of $X$. By Thomason~\cite[Theorem~3.15]{Thomason97} (or by
Balmer~\cite[Theorem~5.5]{Balmer05}) the 
thick tensor ideal in $\dperf X$ generated by $H$ is all
of $\dperf X$. In particular the sheaf $\co_X^{}$ belongs to
the thick tensor ideal generated by $H$; by 
Lemma~\ref{L1.1} there exists an object 
$C\in\dperf X$ and
an integer $L$ so that $\co_X^{}\in\langle C\oo H\rangle_L$. By
Corollary~\ref{C1.17.5} we have $\langle C\oo H\rangle_L\subset
\Coprod_{2L}\big((C\oo H)(-\infty,\infty)\big)$, and hence
\begin{eqnarray*}
\Dqc(X) &=& \Dqc(X)\oo_X^\LL\co_X^{} \\
       &\subset &\Dqc(X)\oo_X^\LL\Coprod_{2L}\big((C\oo_X^\LL H)(-\infty,\infty)\big)\\
   &\subset &\Coprod_{2L}\big(\Dqc(X)\oo_X^\LL C\oo_X^\LL H\big)\\
&\subset &\Coprod_{2L}\big(\Dqc(X)\oo_X^\LL H\big)\\
&\subset &\Coprod_{4L(M+N)}\big((\R f_*G\oplus\R i_*G')(-\infty,\infty)\big)
\end{eqnarray*}
completing the proof that Theorem~\ref{T27.5}
follows from Theorem~\ref{T27.1}.\hfill{$\Box$}

\section{Approximation in the case of quasiprojective schemes}
\label{S2}

We begin by reminding the reader of some standard facts.

\rmd{R2.1}
Let $R$ be a ring. We have an exact sequence $0\la R[x]\stackrel x\la R[x]\la 
R\la 0$. In other words: over the ring $R[x]$ the natural map is a
quasi-isomorphism
\[
\CD
@>>> 0 @>>> R[x] @>x>> R[x] @>>> 0 @>>> \\
@. @VVV @VVV @VVV @VVV @. \\
@>>> 0 @>>> 0       @>>>     R @>>> 0 @>>> 
\endCD
\]
Let us denote the top row 
$\mathrm{Cone}\big(R[x]\stackrel x\sr R[x]\big)$, and the quasi-isomorphism as
$\mathrm{Cone}\big(R[x]\stackrel x\sr R[x]\big)\la R$.
Tensoring this quasi-isomorphism (over $R$) with itself, $n+1$ times, we 
obtain a quasi-isomorphism
\[
\CD
\ds\bigotimes_{i=0}^n\mathrm{Cone}\big(R[x_i]\stackrel {x_i}\sr
R[x_i]\big)
@>>> R
\endCD
\]
To express it slightly differently: over the polynomial ring
$R[x_0^{},x_1^{},\ldots,x_n^{}]$, the Koszul complex on the sequence
$(x_0^{},x_1^{},\ldots,x_n^{})$ is a free resolution of the 
$R[x_0^{},x_1^{},\ldots,x_n^{}]$--module $R$.

Under the standard correspondence, which takes graded
$R[x_0^{},x_1^{},\ldots,x_n^{}]$--modules to
sheaves on $\pp^n_R$, we deduce an exact sequence on $\pp^n_R$
\[\xymatrix@C-4pt{
0\ar[r] & \co(-n) \ar[r] & \co(-n+1)^{\oplus a_{n-1}} \ar[r] &\cdots 
\ar[r]& \co(-1)^{\oplus a_{-1}^{}} \ar[r] & \co^{\oplus a_{0}^{}}
\ar[r] & 
\co(1)\ar[r] & 0
}\]
where the $a_i$ are suitable integers. This gives a quasi-isomorphism
\[\xymatrix@C-4pt{
0\ar[r] & \co(-n) \ar[r]\ar[d] & \co(-n+1)^{\oplus a_{n-1}} \ar[r] \ar[r]\ar[d] &\cdots 
\ar[r]& \co(-1)^{\oplus a_{-1}^{}} \ar[r] \ar[r]\ar[d] & \co^{\oplus a_{0}^{}}
\ar[r] \ar[r]\ar[d] & 0\\
0\ar[r] &0 \ar[r] &0 \ar[r] &\cdots 
\ar[r]& 0 \ar[r] & \co(1)
\ar[r] & 0
}\]
Now it's time to prove something new.
\ermd

\lem{L2.3}
Let $R$ be a ring. In the derived category $\dcoh(\pp^n_R)=
\Dqc(\pp^n_R)^c$ consider the
full subcategory $\cb$ whose set of objects is $\{\co(i),\,-n\leq
i\leq 0\}$. 
Then every line bundle $\co(N)$ belongs to
$\Coprod_{2(n+1)}\big(\cb[-n-1,n]\big)$.
\elem

\prf
The line bundles $\{\co(N),\,-n\leq
N\leq 0\}$ all belong to $\cb\subset
\Coprod_{2(n+1)}\big(\cb[-n-1,n]\big)$, hence we only need to prove
something
for $N\notin[-n,0]$. Let us begin with $N>0$.

In Reminder~\ref{R2.1} we recalled the quasi-isomorphism
\[\xymatrix@C-4pt{
0\ar[r] & \co(-n) \ar[r]\ar[d] & \co(-n+1)^{\oplus a_{n-1}} \ar[r] \ar[r]\ar[d] &\cdots 
\ar[r]& \co(-1)^{\oplus a_{-1}^{}} \ar[r] \ar[r]\ar[d] & \co^{\oplus a_{0}^{}}
\ar[r] \ar[r]\ar[d] & 0\\
0\ar[r] &0 \ar[r] &0 \ar[r] &\cdots 
\ar[r]& 0 \ar[r] & \co(1)
\ar[r] & 0
}\]
Hence $\co(N)=\co(1)^{\otimes N}$ is quasi-isomorphic to the $N\mth$
tensor
power of the quasi-isomorphic complex in the top row.
That is, $\co(N)$ is quasi-isomorphic to a complex 
\[\xymatrix@C-4pt{
\cdots\ar[r] & \co(-n) ^{\oplus b_{n}}  \ar[r] & \co(-n+1)^{\oplus b_{n-1}} \ar[r] &\cdots 
\ar[r]& \co(-1)^{\oplus b_{-1}^{}} \ar[r] & \co^{\oplus b_{0}^{}}
\ar[r] & 
0
}\]
If we take the brutal truncation, that is we define the complex $C$ to
be
\[\xymatrix@C-4pt{
0\ar[r] & \co(-n) ^{\oplus b_{n}}  \ar[r] & \co(-n+1)^{\oplus b_{n-1}} \ar[r] &\cdots 
\ar[r]& \co(-1)^{\oplus b_{-1}^{}} \ar[r] & \co^{\oplus b_{0}^{}}
\ar[r] & 
0
}\]
there are two nonzero homology sheaves: $\ch^0(C)=\co(N)$, while
$\ch^{-n}(C)=\ce$
is some vector bundle on $\pp^n_R$. The \tsrt\ truncation gives a
triangle
$\T^{n}\ce\la C\la \co(N)\la\T^{n+1}\ce$, and the morphism
$\co(N)\la\T^{n+1}\ce$ belongs to
$\Ext^{n+1}_{\pp^n_R}(\co(N),\ce)=H^{n+1}\big(\ce(-N)\big)=0$. 
Therefore $C\cong\co(N)\oplus\T^n\ce$.
By Example~\ref{E1.18} the object $C$ belongs to
$\Coprod_{n+1}\big(\cb[-n,0]\big)$,
and 
the line bundle $\co(N)$, being a direct summand of $C$, belongs to
$\Smr\Big[\Coprod_{n+1}\big(\cb[-n,0]\big)\Big]\subset
\Coprod_{2(n+1)}\big(\cb[-n-1,0]\big)$. The last inclusion is by
Lemma~\ref{L1.15}. 

It remains to consider the case $N<-n$. In the above we produced, for
every
$M>0$, a
complex $C$ of the form
\[\xymatrix@C-4pt{
0\ar[r] & \co(-n) ^{\oplus b_{n}}  \ar[r] & \co(-n+1)^{\oplus b_{n-1}} \ar[r] &\cdots 
\ar[r]& \co(-1)^{\oplus b_{-1}^{}} \ar[r] & \co^{\oplus b_{0}^{}}
\ar[r] & 
0
}\]
which is quasi-isomorphic to $\co(M)\oplus \T^n\ce$ where $\ce$ is a vector
bundle.
Then $\RHHom\big(C,\co(-n)\big)$ is quasi-isomorphic to
$\co(-n-M)\oplus\T^{-n}\RHHom\big(\ce,\co(-n)\big)$.
This makes $\co(-n-M)$ a direct summand of
$\RHHom\big(C,\co(-n)\big)\in
\Coprod_{n+1}\big(\cb[0,n]\big)$. Any $N<-n$ can be written
as $N=-n-M$ for $M>0$, and the 
line bundle $\co(N)=\co(-n-M)$ belongs to
$\Smr\Big[\Coprod_{n+1}\big(\cb[0,n]\big)\Big]\subset
\Coprod_{2(n+1)}\big(\cb[-1,n]\big)$.
\eprf

\cor{C2.5}
Let $R$ be a ring, $X$ a scheme and $f:X\la\pp^n_R$ a morphism
(a base-point free
map). Let $\cl=f^*\co(1)$.
In the category $\Dqc(X)$ define $\cb$ to be the full subcategory with
the set of
objects $\{\cl^{i},\,-n\leq i\leq 0\}$.

Then every line bundle $\cl^N$ belongs to $\Coprod_{2(n+1)}\big(\cb[-n-1,n]\big)$.
\ecor

\prf
Pull back the inclusion of Lemma~\ref{L2.3} via $\LL f^*$.
\eprf

Now we are ready to prove Theorem~\ref{T0.3} for quasiprojective $X$.
Note that the bounds are quite explicit and simple.

\pro{P2.7}
Let the notation be as in Corollary~\ref{C2.5}, but assume furthermore
that $\cl$ is ample on $X$. Suppose $F$ is an
object of $\Dqc(X)^{\leq B}$, that
is $\ch^i(F)=0$ for all $i>B$,  and let $A\leq B$ be an integer.

Then there exists a triangle $D\la E\la F\la$ in $\Dqc(X)$,
with $D\in\Dqc(X)^{<A}$  and
$E\in\Coprod_{(B-A+2)(2n+2)}\big(\cb[A-n-2,B+n]\big)$. If $X$ is noetherian
and the cohomology of $F$ is coherent, then $E$ may be chosen to be
compact.
\epro

\prf
Let $\cc$ be the full subcategory of the category of line bundles on $X$, whose
set of objects is $\{\cl^i,\,i\in\zz\}$. Because $\cl$ is ample,
Illusie~\cite[Lemme~2.2.8]{Illusie71B} tells us that
the complex $F$ has a resolution in terms of the $\cl^i$: there exists
a cochain complex
\[\xymatrix{
\cdots\ar[r]&  
x^{B-2}_{}\ar[r] &
x^{B-1}_{}\ar[r] & x^{B}_{}\ar[r] &0\ar[r] &\cdots
}\]
quasi-isomorphic to $F$, and where the $x^i$ are all coproducts of
objects in $\cc$. If $X$ is noetherian and $F$ has coherent cohomology 
then Illusie~\cite[Proposition~2.2.9]{Illusie71B} allows us to choose
all the $x^i$ to be finite coproducts of objects of $\cc$. Let $E$ be 
the brutal truncation 
\[\xymatrix{
\cdots\ar[r]& 0 \ar[r] &x^{A-1}_{}\ar[r] &
x^{A}_{}\ar[r] &\cdots\ar[r] &
x^{B-1}_{}\ar[r] & x^{B}_{}\ar[r] &0\ar[r] &\cdots
}\]
In the triangle $D\la E\la F\la$ we have $D\in\Dqc(X)^{<A}$, and in the 
case where $X$ is noetherian and 
the cohomology of $F$ is coherent the object $E$ is compact. 
Example~\ref{E1.18} tells us that $E\in\Coprod_{B-A+2}\big(\cc[A-1,B]\big)$,
while Corollary~\ref{C2.5} guarantees that
$\cc[A-1,B]\subset\Coprod_{2(n+1)}\big(\cb[A-n-2,B+n]\big)$. The 
Proposition follows.
\eprf

\section{Approximation in the general case}
\label{S3}

\rmd{R3.1}
We recall \cite[Theorem~4.1 and Theorem~4.2]{Lipman-Neeman07}; 
we only care here about the
case where $X$ is noetherian and separated. 
\cite[Theorem~4.2]{Lipman-Neeman07} tells us that the category 
$\Dqc(X)$ has a compact generator $G$, while 
\cite[Theorem~4.1]{Lipman-Neeman07} asserts that every object 
in $\dmcoh(X)$ can be approximated by compacts arbitrarily well. 
This means: for any object $F\in\dmcoh(X)$ and any integer
$m$, there exists a triangle $D\la E\la F\la$ with
$E$ a compact object of
$\Dqc(X)$ and $D\in\dmcoh(X)^{<m}$. If we abuse the notation 
slightly and write $G$ for the category 
with one object $G$, \cite[Theorem~4.2]{Lipman-Neeman07} tells us 
that the object $E$ of \cite[Theorem~4.1]{Lipman-Neeman07} lies in 
\[
\Smr\Big[\COprod\big(G(-\infty,\infty)\big)\Big]=\Smr\left[\bigcup_{n>0,A\leq B}
\COprod_n\big(G[A,B]\big)\right].
\]
It follows that, for some integers $n>0$ and $A\leq B$, the 
object $E$ must belong to $\Smr\left[
\COprod_n\big(G[A,B]\big)\right]\subset\Coprod_{2n}\big(G[A-1,B]\big)$
where the inclusion is by Lemma~\ref{L1.15}.
\ermd

This motivates the following

\dfn{D3.3}
Let $X$ be a noetherian, separated scheme.
We say that a subcategory $\cs\subset\Dqcmi(X)$ is \emph{approximable} 
if there exists a compact generator $G\in\Dqc(X)$ with the
property that,
for every integer $m$, there exist integers
$n,A\leq B$, depending only on $\cs$, $G$ and $m$, so that any
object $F\in\cs$ admits a triangle $D\la E\la F\la$ with $D\in\Dqc(X)^{<m}$
and $E\in\Coprod_n\big(G[A,B]\big)$.
\edfn

\rmk{R3.4}
If $\cs\subset\Dqcmi(X)$ 
is approximable for one compact generator it is approximable for every 
compact generator. If $G,G'$ are two compact generators, the fact that 
$G'$ is compact while $G$ generates tells us that for some
integers $n,A,B$ the object $G'$ belongs to
$\Smr\left[
\COprod_n\big(G[A,B]\big)\right]\subset\Coprod_{2n}\big(G[A-1,B]\big)$,
and by symmetry $G$ belongs to
$\Smr\left[
\COprod_{n'}\big(G'[A',B']\big)\right]\subset\Coprod_{2n'}\big(G'[A'-1,B']\big)$.
Any $G$--approximation $D\la E\la F\la $ of an object in $F\in\cs$ 
is also a 
$G'$--approximation, only the integers change.
\ermk

\rmk{E3.4.4}
Let $G$ be a classical generator of $\D^{\text{perf}}(X)$ [see
  Reminder~\ref{R0.-3} for the definition] and let $G'$ be
a compact generator in $\Dqc(X)$. Because $G'$ is compact
it lies in $\D^{\text{perf}}(X)$, and as $G$ is a 
classical generator of $\D^{\text{perf}}(X)$
we have $G'\in \Smr\left[
\COprod_n\big(G[A,B]\big)\right]$ for some integers $n,A,B$.
Hence $G$ must be a compact generator for $\Dqc(X)$.
The statements in the introduction assert that classical generators
of $\D^{\text{perf}}(X)$ have some properties, and it suffices to show 
that compact generators of $\Dqc(X)$ have these properties.
\ermk

\exm{E3.4.5}
Let $X$ be a scheme over a ring $R$, with an ample line bundle $\cl$. 
Suppose that some sections of $\cl$ give a morphism 
(i.e. a base-point free map)
$f:X\la\pp^n_R$, and let $G=\oplus_{i=-n}^0\cl^i$. Clearly $G$ is compact.

Since $\cl^i,\,-n\leq i\leq 0$ are direct summands of $G$ they belong to
$\Smr\left[
\COprod_1\big(G[0,0]\big)\right]\subset\Coprod_{2}\big(G[-1,0]\big)$.
In the notation of Corollary~\ref{C2.5} we have that
$\cb\subset\Coprod_{2}\big(G[-1,0]\big)$, and the Corollary
tells us that every line bundle of the form
$\cl^N$ belongs to
\[
\Coprod_{2(n+1)}\big(\cb[-n-1,n]\big)\sub
\Coprod_{4(n+1)}\big(G[-n-2,n]\big)\ .
\]
By \cite[Example~1.10]{Neeman96} or~\cite[Tag~0BQQ]{stacks-project}
the line bundles $\cl^N$ weakly generate $\Dqc(X)$, and from 
\cite[Lemma~3.2]{Neeman96} it follows
that they generate. Since the $\cl^n$ all belong to the 
subcategory generated by $G$, we conclude that $G$ is a 
compact generator. 

Proposition~\ref{P2.7} now tells us that the category
$\Dqc(X)^{\leq B}$ is approximable for every integer $B$.
\eexm

\lem{L3.5}
If $\cs,\cs'$ are approximable subcategories of $\Dqcmi(X)$, then so is
$\cs\star\cs'$.
\elem

\prf
Let $m$ be an integer and $G$ a compact generator for
$\Dqc(X)$. Because $\cs'$ is approximable there are
 integers
$n',A',B'$ so that any $F'\in\cs'$ admits a triangle
$D'\la E'\la F'\la$ with $D'\in\Dqc(X)^{<m}$ and
$E'\in\Coprod_{n'}\big(G[A',B']\big)$.

Next by Hartshorne~\cite[Exercise~III.4.8(c)]{Hartshorne77}
we may
choose an integer $\ell$ such that $H^i(X,\ca)=0$ for all quasicoherent
sheaves $\ca$ and all $i>\ell$.
Now $G$ is compact, meaning a perfect complex, and hence 
$G^\vee=\RHHom(G,\co_X^{})$ is also a perfect complex, and belongs to
$\Dqc(X)^{\leq a}$ for some integer $a$. If $L$ belongs to $\Dqc(X)^{\leq b}$
then $\RHHom(G,L)=G^\vee\oo L\in\Dqc(X)^{\leq a+b}$, and 
$\Hom(G,L)=H^0\big[\RHHom(G,L)\big]=0$ whenever $a+b<-\ell$. That is
$\Hom\Big[G,\Dqc(X)^{\leq-a-\ell-1}\big]=0$, and hence
\[
\Hom\Big[\Coprod_{n'}\big(G[A',B']\big)\,\,,\,\,
\Dqc(X)^{\leq A'-a-\ell-1}\Big]=0\ .
\]

Now let $\wt m=\min(m,A'-a-\ell)$. Because $\cs$ is approximable 
there exist integers $n,A,B$ so that any object $F\in\cs$ admits
a triangle $D\la E\la F\la $ with $E\in \Coprod_{n}\big(G[A,B]\big)$
and $D\in\Dqc(X)^{<\wt m}$.
So much for the construction: I assert that every object $\wt F\in\cs\star\cs'$
admits a triangle $\wt D\la \wt E\la\wt F\la$, with 
$\wt D\in\Dqc(X)^{<m}$ and 
$\wt E\in\Coprod_{n+n'}\Big(G\big[\min(A,A'),\max(B,B')\big]\Big)$.
It remains to prove the assertion.

Let $\wt F$ be an object of $\cs\star\cs'$. Then there exists a triangle
$F\la \wt F\la F'\la$ with $F\in\cs$ and $F'\in\cs'$. By the above there 
exist triangles $D\la E\la F\la $ and $D'\la E'\la F'\la $ with
$E\in\Coprod_{n}\big(G[A,B]\big)$, $E'\in\Coprod_{n'}\big(G[A',B']\big)$,
$D\in\Dqc(X)^{<\wt m}$ and $D'\in\Dqc(X)^{< m}$. We have a diagram
\[\xymatrix{
E' \ar[r] & F' \ar[r]\ar[d] & \T D'\\
\T E \ar[r] &\T F \ar[r]& \T^2 D
}\]
where the rows are triangles, 
and the composite from top left to bottom right is a map
from $E'\in\Coprod_{n'}\big(G[A',B']\big)$ to
$\T^2D\in\Dqc(X)^{<\wt m-2}\subset \Dqc(X)^{<A'-a-\ell-2}$. Since
\[
\Hom\Big[\Coprod_{n'}\big(G[A',B']\big)\,\,,\,\,
\Dqc(X)^{\leq A'-a-\ell-1}\Big]=0
\]
the map $E'\la \T^2D$ must vanish, and there is a map $E'\la\T E$ rendering 
commutative the square
\[\xymatrix{
E' \ar[r]\ar[d] & F' \ar[d] \\
\T E \ar[r] &\T F 
}\]
Now complete this commutative square to a $3\times3$ diagram: we obtain a 
diagram where the rows and columns are triangles
\[\xymatrix{
D \ar[r]\ar[d] & E\ar[r]\ar[d] & F \ar[r]\ar[d]& \\
\wt D \ar[r]\ar[d] &\wt E\ar[r]\ar[d] &\wt F \ar[r]\ar[d]& \\
D' \ar[r]\ar[d] & E'\ar[r]\ar[d] & F' \ar[r]\ar[d]& \\
 & & & 
}\]
 From the triangle $E\la \wt E \la E'\la$ and the fact that $E\in\Coprod_{n}\big(G[A,B]\big)$ and
$E'\in\Coprod_{n'}\big(G[A',B']\big)$ we learn that 
$\wt E\in\Coprod_{n+n'}\Big(G\big[\min(A,A'),\max(B,B')\big]\Big)$.
The triangle $D\la\wt D\la D'\la$ and the 
fact that $D\in \Dqc(X)^{<\wt m}\subset\Dqc(X)^{< m}$ and 
$D'\in\Dqc(X)^{< m}$ give
$\wt D\in \Dqc(X)^{< m}$. And the triangle $\wt D\la\wt E\la\wt F\la$ 
is now as desired.
\eprf

\lem{L3.7}
$f:X\la Y$ be a proper map of noetherian schemes, and 
let $\cs$ be an approximable subcategory of $\Dqcmi(X)$. 
Then $\R f_*\cs$ is an approximable
subcategory of $\Dqcmi(Y)$.
\elem

\prf
Choose a compact generator $G\in\Dqc(X)$. Because $f$ is proper
the object $\R f_*G$ has bounded-above coherent cohomology; it 
belongs to $\dmcoh(Y)$ and 
in Reminder~\ref{R3.1} we learned that it is approximable. 
Slightly more generally: given any pair of integers
$A\leq B$,
the finite category with objects $\cb=\{\T^i\R f_*G,\,-B\leq i\leq -A\}$
is approximable. By Lemma~\ref{L3.5} so is
$\Coprod_n(\cb)$ for any $n$.

Choose an integer $\ell$ so that $\R f_*\Dqc(X)^{\leq0}\subset
\Dqc(Y)^{\leq\ell}$. Given an integer $m$, the fact that 
$\cs$ is approximable in $\Dqcmi(X)$ allows us to choose 
integers $n,A,B$ so that
any
object $F\in\cs$ admits a triangle $D\la E\la F\la$ with $D\in\Dqc(X)^{<m-\ell}$
and $E\in\Coprod_n\big(G[A,B]\big)$. 
By the first paragraph the category $\Coprod_n\big(\R f_*G[A,B]\big)$
is approximable in $\Dqcmi(Y)$; fix a compact generator 
$H\in\Dqc(Y)$, and we may choose integers $n',A',B'$ so that
every object $F'\in \Coprod_n\big(\R f_*G[A,B]\big)$
admits a triangle $D'\la E'\la F'\la$ with $D'\in\Dqc(Y)^{<m}$
and $E'\in\Coprod_{n'}\big(H[A',B']\big)$.
We have now chosen all our integers. I assert that every 
object $\R f_*F\in\R f_*\cs$ admits a triangle
$\wt D\la\wt E\la \R f_*F\la$, with
$\wt D\in \Dqc(Y)^{<m}$ and $\wt E\in\Coprod_{n'}\big(H[A',B']\big)$.

It remains to prove the assertion. By our choice of
$n,A,B$, the object  $F\in\cs$ admits a triangle 
$D\la E\la F\la$ with $D\in\Dqc(X)^{<m-\ell}$
and $E\in\Coprod_n\big(G[A,B]\big)$.
In $\Dqc(Y)$ we deduce a triangle
$\R f_*D\la\R f_* E\la \R f_*F\la$.
The choice of the integer $\ell$ guarantees that
 $\R f_*D\in \Dqc(Y)^{<m}$, and we know that
\[
\R f_*E\in\R f_*\Coprod_n\big(G[A,B]\big)\subset
\Coprod_n\big(\R f_*G[A,B]\big)\ .
\]
By our choice of the integers $n',A',B'$ there exists a triangle
$D'\la\wt E\la \R f_* E\la $ with
$D'\in\Dqc(Y)^{<m}$ and $\wt E\in\Coprod_{n'}\big(H[A',B']\big)$.
Now complete the composable maps $\wt E\la\R f_*E\la\R f_*F$ to
an octahedron; we obtain a diagram where the rows and columns are
triangles
\[\xymatrix{
D' \ar[r]\ar@{=}[d] & \wt D \ar[r]\ar[d] & \R f_*D\ar[r]\ar[d] & \\
D' \ar[r] & \wt E \ar[r]\ar[d] & \R f_*E\ar[r]\ar[d] & \\
& \R f_*F \ar@{=}[r]\ar[d] & \R f_*F\ar[d] &\\
 & &  &
}\]
The fact that $D',\R f_*D$ both belong to $\Dqc(Y)^{<m}$ says
that so does $\wt D$, and the triangle $\wt D\la\wt E\la\R f_*F\la$ has
the desired properties.
\eprf

And now for the proof of Theorem~\ref{T0.3}. The statement, rephrased in terms
of the categories $\Coprod_N\big(G[A,B]\big)$ as in Remark~\ref{R1.17.7}, is

\thm{T3.9}
Let $X$ be a separated scheme, of finite type over a
noetherian ring $R$. Then the category $\Dqc(X)^{\leq n}$ is approximable
for any integer $n$.
\ethm

\prf
If the theorem were false, then the set of closed subschemes of $X$ on which
it fails would have a smallest member. We may therefore assume the
theorem is true on every proper, closed subscheme of $X$.

By Chow's Lemma we may
choose a proper morphism $f:Y\la X$, where $Y$ is quasiprojective over $R$
and $f$ is an isomorphism on a dense open subset. In the triangle
$Q\la \co_X\la\R f_*\LL f^*\co_X\la $ we have that $Q$ is isomorphic
to a bounded complex of coherent sheaves, which vanish on the dense
open set where $f$ is an isomorphism. By~\cite[Lemma~7.40]{Rouquier08} 
there is a closed immersion
$i:Z\la X$ and an object of $P\in\dcoh(Z)$, so that
$Q=\R i_*P$.

Now take any object in $F\in\Dqc(X)$; tensoring the triangle with $F$
we obtain a triangle $F\oo_X^\LL \R i_*P\la F \la F\oo_X^\LL \R f_*\co_Y^{}\la$,
and the projection formula allows us to rewrite this as
$\R i_*[\LL i^*F\oo P]\la F\la \R f_*\LL f^* F\la$. If $F\in\Dqc(X)^{<m}$ then
$\LL f^*F\in\Dqc(Y)^{<m}$ and $\LL i^*F\in\Dqc(Z)^{<m}$, and since 
$P$ is bounded it lies in some $\Dqc(Z)^{\leq\ell}$ and hence
$\LL i^*F\oo P\in\Dqc(Z)^{<m+\ell}$. The triangle therefore shows
that $\Dqc(X)^{<m}\subset\big[\R i_*\Dqc(Z)^{<m+\ell}\big]\star
\big[\R f_*\Dqc(Y)^{<m}\big]$.

Because $Z$ is a proper closed subscheme of $X$ the category $\Dqc(Z)^{<m+\ell}$
is approximable in $\Dqcmi(Z)$, while the category $\Dqc(Y)^{<m}$ is
approximable in $\Dqcmi(Y)$ by Example~\ref{E3.4.5}. Lemma~\ref{L3.7} tells
us that $\R i_*\Dqc(Z)^{<m+\ell}$ and $\R f_*\Dqc(Y)^{<m}$ are approximable,
and Lemma~\ref{L3.5} permits us to conclude that $\Dqc(X)^{<m}$ is
approximable.
\eprf

\section{Proofs of Theorems~\protect{\ref{T0.11}} and~\protect{\ref{T27.1}}}
\label{S4}

We begin with a straightforward corollary of Theorem~\ref{T3.9}: 

\cor{C0.7}
Let $m\leq n$ be integers, let $X$ be a scheme of finite type over
a noetherian ring $R$, and let $G$ be a classical generator for
$\dperf X\subset\Dqc(X)$. There exist integers $N,A,B$,
depending only on $G$, $m$ and $n$, so that if $F\in\Dqc(X)$ satisfies
\be
\item
$F\in\Dqc(X)^{\leq n}$, and 
\item
$\RHHom(F,Q)\in
\Dqc(X)^{\leq -m}$ for every $Q\in\Dqc(X)^{\leq0}$,
\ee
then $F$ belongs to
$\Coprod_N\big(G[A,B]\big)$.
\ecor

\prf
Choose an integer $\ell$ so that $H^i(X;\cm)=0$ for all quasicoherent
sheaves $\cm$ and all $i>\ell$, and then apply Theorem~\ref{T3.9} to
the integers $\min(n,m-\ell)\leq n$, the scheme $X$ and the
classical generator 
$G\in \dperf X$. There exist integers $N,A\leq B$ so that any
$F\in\Dqc(X)^{\leq n}$ admits a triangle
$D\la E\la F\la $ with $E$ in $\Coprod_N\big(G[A,B]\big)$
and
$D\in \Dqc(X)^{< m-\ell}$. 
Now assume that $F$ also satisfies 
condition (ii) of the Corollary.
The map $F\la \T D$ is an element in
\begin{eqnarray*}
\Hom(F,\T D) &\subset& H^0\big[\RHHom(F,\T D )\big]\\
&\subset& H^0\big[ \Dqc(X)^{<-\ell}\big]\\
&=& 0\ .
\end{eqnarray*}
Hence the map $F\la\T D$ vanishes, and $F$ is a direct summand of
$E\in \Coprod_N\big(G[A,B]\big)$. Lemma~\ref{L1.15} tells
us that $F$ belongs to
$\Coprod_{2N}\big(G[A-1,B]\big)$.
\eprf

In Sections~\ref{S2} and \ref{S3} our schemes were mostly assumed noetherian.
The next result drops this hypothesis: we're about to prove
Theorem~\ref{T0.11}. We state its
$\Coprod_N\big(G[A,B]\big)$ version for the reader's convenience.
In the proof we will freely appeal
to the fact that the categories $\Dqc(X)$ 
contain compact generators for any quasicompact, quasiseparated scheme
$X$---see \cite[Theorem~3.1.1]{BondalvandenBergh04} or
\cite[Theorem~4.2]{Lipman-Neeman07}.

\thm{T4.1}
Let $j:V\la X$ be an open immersion of quasicompact, separated schemes,
and let
$G$ be a compact generator for $\Dqc(X)$.
If $H$ is any compact object of $\Dqc(V)$, and we are given
integers $n,a\leq b$, then there exist integers
$N,A\leq B$ so that
$\Coprod_n\big(\R j_*H[a,b]\big)\subset\Coprod_N\big(G[A,B]\big)$.
\ethm

\prf
Since $X$ is quasicompact we may write it as a
finite union $X=\cup_{i=1}^\ell U_i$
of open affines $U_i$. Let $X_k=V\cup\big[\cup_{i=1}^kU_i\big]$, then
we may factor $j$ as the composite of the inclusions
\[\xymatrix{
  U=X_0\ar[r]^-{j_i^{}} &X_1\ar[r]^-{j_2^{}} &X_2\ar[r]^-{j_3^{}} &\cdots
  \ar[r]^-{j_{\ell-1}^{}} & X_{\ell-1}\ar[r]^-{j_\ell^{}} & X_\ell=X
}\]
and it suffices to prove the theorem for each $j_i^{}:X_{i-1}\la X_i$.
Thus we may
suppose $X=U\cup V$
expresses $X$ as the union of the quasicompact open subsets $U$ and $V$,
with $U$ affine. 
By Thomason's localization theorem we may choose a compact object 
$\wt H\in\Dqc(X)$ and a quasiisomorphism $\LL j^*\wt H\cong H\oplus \T H$;
see Thomason and Trobaugh~\cite{ThomTro} or 
\cite[statements~2.1.4 and 2.1.5]{Neeman96}. Now consider the triangle
\[
\CD
Q @>>> \co_X^{} @>>> \R j_*\LL j^*\co_X^{} @>>> \T Q
\endCD
\]
If we tensor with $\wt H$, and use the projection formula to obtain
the first isomorphism in 
$[\R j_*\LL j^*\co_X^{}]\oo_X^\LL\wt H\cong \R j_*\LL j^*\wt H
\cong \R j_*(H\oplus\T H)$,
 we obtain a triangle
\[
\CD
Q\oo^\LL_X\wt H @>>> \wt H @>>> \R j_*H\oplus\T\R j_*H @>>> \T Q\oo^\LL_X\wt H
\endCD
\]
Because $\wt H$ is compact in $\Dqc(X)$ and $G$ is a compact generator
there exist integers $N,A\leq B$ with $\wt H\in\Coprod_N\big(G[A,B]\big)$.
 From the triangle it suffices to show that there also exist integers 
$N',A'\leq B'$ with $Q\oo^\LL_X\wt H\in \Coprod_{N'}\big(G[A',B']\big)$.

Now $Q$ vanishes on the open set $V$; if $h:U\la X$ is the inclusion
then the natural map $Q\la \R h_*\LL h^*Q$ is an isomorphism. This makes
$Q\oo^\LL_X\wt H\cong [\R h_*\LL h^*Q]\oo^\LL_X\wt H\cong
Q\oo^\LL_X[\R h_*\LL h^*\wt H]$. But $\LL h^*\wt H$ is a compact
object in $\Dqc(U)$ and $\co_U^{}$ is a compact generator (this is because
$U$ is affine), and hence $\LL h^*\wt H$ must belong to $\Coprod_M\big(
\co_U^{}[A'',B'']\big)$ for some integers $M,A''\leq B''$. 
Applying the functor $Q\oo^\LL_X\R h_*(-)$ tells us that
$Q\oo^\LL_X[\R h_*\LL h^*\wt H]$ belongs
to  $\Coprod_M\big(
[Q\oo^\LL_X\R h_*\co_U^{}][A'',B'']\big)$,
 and it
suffices to show that $Q\oo^\LL_X\R h_*\co_U^{}$ belongs to
$\Coprod_{N}\big(G[A,B]\big)$ for some $N,A\leq B$. Now observe 
the isomorphisms
 $Q\oo^\LL_X\R h_*\co_U^{}\cong Q\oo^\LL_X\R h_*\LL h^*\co_X^{}\cong
\R h_*\LL h^*Q\cong Q$; it remains to show that
$Q$ belongs to  $\Coprod_N\big(G[A,B]\big)$ for some $N,A\leq B$.

Now $X$ is separated and $U,V$ are quasicompact open subsets of $X$; hence
$U\cap V$ is quasicompact. It is a quasicompact open subset of the open
affine subset $U=\spec S$. Absolute noetherian approximation, that is 
Thomason and Trobaugh~\cite[Theorem~C.9]{ThomTro} or~\cite[Tags 01YT
and 081A]{stacks-project},
allows us to choose a scheme $Y$ of finite type over $\zz$ and an affine
map $f:X\la Y$, so that
\be
\item
$Y$ is separated.
\item
There is an open affine subset $U'\subset Y$ and an open subset $V'\subset
Y$ with $U'\cup V'=Y$, and
so that $f^{-1}U'=U$ and $f^{-1}V'=V$. Let $h':U'\la Y$, $j':V'\la Y$
be the open immersions.
\ee

The affine scheme $U'$ may be written as $U'=\spec{S'}$, with $S'$ a
noetherian ring. The set
$Y-V'=U'-U'\cap V'$ is a Zariski closed subset of the noetherian
affine
scheme $U'$, hence there are elements
$\{g'_1,g'_2,\ldots,g'_n\}\subset S'$ so that $Y-V'$ is precisely
the subset $V(g'_1,g'_2,\ldots,g'_n)$ of $U'$ on which the $g'_i$ all vanish. 
In the category $\D(S')\cong\Dqc(U')$
consider the complexes $L_i=\{0\la S'\la S'[{g'_i}^{-1}]\la 0\}$, with $S'$ in
degree 0 and $S'[{g'_i}^{-1}]$ in degree 1. Put
$Q'=\otimes_{i=1}^nL_i$; the natural cochain map $Q'\la S'$ may be
completed 
to a triangle $Q'\la S'\la T\la$ in $\D(S')$, which we recognize in 
$\D(S')\cong\Dqc(U')$ as the canonical triangle
\[
\CD
Q' @>>> \co_{U'}^{} @>>> \R\alpha_*\co_{U'\cap V'}^{} @>>> 
\endCD
\]
where $\alpha:U'\cap V'\la U'$ is the open immersion. Since $\R h'_*Q'$
is the extension by zero of $Q'$ to all of $Y$
 we obtain a triangle in $\Dqc(Y)$
\[
\CD
\R h'_*Q' @>>> \co_{Y}^{} @>>> \R j'_*\co_{ V'}^{} @>>> 
\endCD
\]
and pulling back via $f:X\la Y$ we deduce that $Q\cong\LL f^*\R h'_*Q'$. 
Note that $Q'$ is the homotopy colimit
over $r$
of the complexes
\[
\oo_{i=1}^n\{S'\stackrel{{g'_i}^r}\la S'\}
\]
and hence fits in a triangle
\[
\CD
\ds\coprod_{r=1}^\infty\left[\oo_{i=1}^n\{S'\stackrel{{g'_i}^r}\la S'\}\right]
@>>> Q'
@>>> \ds\T \coprod_{r=1}^\infty\left[\oo_{i=1}^n\{S'\stackrel{{g'_i}^r}\la S'\}\right]@>>>
\endCD
\]
This means that the object $Q'\in\D(S')$ has a projective 
resolution that vanishes outside the interval $[-1,n]$. Let $\ell$ be
the open immersion $\ell:U'\cap V'\la Y$.
If $F$ is any
 object in $\Dqc(Y)^{\leq0}$, we have the standard triangle 
\[
\CD
F @>>> \R h'_*\LL {h'}^* F\oplus \R j'_*\LL {j'}^*F @>>>\R\ell_*\LL\ell^*F @>>> 
\endCD
\]
and, applying the functor $\RHHom(\R h'_*Q',-)$, 
we obtain a triangle
\[
\CD
\RHHom(\R h'_*Q',F) @>>> \RHHom(\R h'_*Q',\R h'_*\LL {h'}^* F)\oplus \RHHom(\R h'_*Q',\R j'_*\LL {j'}^*F) \\
@. @VVV \\
@. \RHHom(\R h'_*Q',\R\ell_*\LL\ell^*F) 
\endCD
\]
The vanishing of the objects
$
\RHHom(\R h'_*Q',\R j'_*\LL {j'}^*F)\cong \R j'_*\RHHom(\LL {j'}^*\R h'_*Q',\LL {j'}^*F)
$
and $\RHHom(\R h'_*Q',\R\ell_*\LL\ell^*F)\cong\R\ell_*
 \RHHom(\LL\ell^*\R h'_*Q',\LL\ell^*F)$ is because $\LL {j'}^*\R h'_*Q'=0=
\LL\ell^*\R h'_*Q'$. From the triangle we learn that $\RHHom(\R h'_*Q',F)$
is isomorphic to $\RHHom(\R h'_*Q',\R h'_*\LL {h'}^* F)=
\R h'_*\RHHom(\LL {h'}^*\R h'_*Q',\LL {h'}^* F)=\R h'_*\RHHom(Q',\LL {h'}^* F)$, which 
we can compute using the projective resolution for $Q'$ on the affine open
set $U'$. We deduce that $\RHHom(\R h'_*Q',F)\in\Dqc(Y)^{\leq1}$ for
every object $F$ in $\Dqc(Y)^{\leq0}$.

Corollary~\ref{C0.7} therefore applies. On the 
separated scheme
$Y$, of finite type over $\zz$, we have that $\R h'_*Q'$ lies
in $\Coprod_K\big(G'[\ov A,\ov B]\big)$ for some integers $K,\ov A\leq\ov B$
and some compact generator $G'$.
Thus $Q=\LL f^*\R h'_*Q'$ lies 
in $\LL f^*\Coprod_K\big(G'[\ov A,\ov B]\big)\subset
\Coprod_K\big(\LL f^*G'[\ov A,\ov B]\big)$.
Since $\LL f^*G'$ is compact and $G$ is a compact generator of $\Dqc(X)$,
we have
$\LL f^*G'\in \Coprod_{K'}\big(G[\alpha,\beta]\big)$ for some
integers $K',\alpha\leq\beta$, and therefore
$Q\in \Coprod_{KK'}\big(G(\alpha+\ov A,\beta+\ov B]\big)$.
\eprf

We finish the article with

\medskip

\nin
\emph{Proof of Theorem~\ref{T27.1}.}\ \ 
If $X$ is affine, that is $X=\spec R$ for some ring $R$ of
finite global dimension, then 
at some level the result goes back to Kelly~\cite{Kelly65}; see also
Street~\cite{Street73}.
The reader can find 
modern treatments in Christensen~\cite[Corollary~8.4]{Christensen96}
or Rouquier~\cite[Proposition~7.25]{Rouquier08}. More precisely: 
if the ring $R$ is of global dimension $\leq M$, it
is classical
 that $\Dqc\big(\spec R\big)\subset \Coprod_{M+2}\big(\co_{\spec R}^{}(-\infty,\infty)\big)$.

We treat the general case by induction on the number 
of open affine subsets in the cover of $X$. 
Suppose we know the theorem for all schemes which admit a cover
by $\leq n$ open affines $U_i=\spec{R_i}$, with each $R_i$ of finite
global dimension. Let $X$ be a scheme admitting a cover by $n+1$ open
affines $U_i=\spec{R_i}$, with each $R_i$ of finite
global dimension. Put $V=\bigcup_{i=1}^{n}U_i$ and 
$U=U_{n+1}$, and let $G_1$ be a compact generator for
$\Dqc(U)$, $G_2$ a compact generator for $\Dqc(V)$ and $G$ a compact
generator for $\Dqc(X)$. Let $j_1:U\la X$ and $j_2:V\la X$ be
the inclusions. By
Theorem~\ref{T4.1} there exists an integer
$N$ so that $\R j_{1*}G_1$ and $\R j_{2*}G_2$ both belong to
$\Coprod_{N}\big(G(-\infty,\infty)\big)$. 
By induction there is an integer $M$ so that 
$\Dqc(U)=\Coprod_{M}\big(G_1(-\infty,\infty)\big)$ and
$\Dqc(V)=\Coprod_{M}\big(G_2(-\infty,\infty)\big)$. Therefore
$\R j_{1*}\Dqc(U)$ and $\R j_{2*}\Dqc(V)$ are both contained in
$\Coprod_{MN}\big(G(-\infty,\infty)\big)$.

But now $X=U\cup V$. Put $W=U\cap V$ and let $j:W\la U$ be the
inclusion. 
Any object $F\in\Dqc(X)$ fits in a triangle
\[
\CD
\R j_{1*}\big[\R j_*\LL j^*\LL j_1^*\Tm F\big]@>>>F 
@>>> \R j_{1*}\big[\LL j_{1}^*F\big]\oplus  \R j_{2*}\big[\LL j_{2}^*F\big]
@>>> 
\endCD
\]
Thus $F$ belongs to $\big[\R j_{1*}\Dqc(U)\big]\star
\big[
\R j_{1*}\Dqc(U)\oplus \R j_{2*}\Dqc(V)
\big]$,
which is contained in
\[
\Coprod_{MN}\big(G(-\infty,\infty)\big)\star
\Coprod_{MN}\big(G(-\infty,\infty)\big)\eq
\Coprod_{2MN}\big(G(-\infty,\infty)\big)\ .
\]
\hfill{$\Box$}

\def\cprime{$'$}
\providecommand{\bysame}{\leavevmode\hbox to3em{\hrulefill}\thinspace}
\providecommand{\MR}{\relax\ifhmode\unskip\space\fi MR }
\providecommand{\MRhref}[2]{%
  \href{http://www.ams.org/mathscinet-getitem?mr=#1}{#2}
}
\providecommand{\href}[2]{#2}

\end{document}